\documentclass[a4paper,11pt]{article}
\usepackage[english]{babel}
\usepackage{amssymb}
\usepackage{amsfonts}
\usepackage{amsmath}

\newtheorem{theorem}{Theorem}

\newtheorem{definition}[theorem]{Definition}

\newtheorem{lemma}[theorem]{Lemma}

\newtheorem{proposition}[theorem]{Proposition}
\newtheorem{remark}[theorem]{Remark}

\input{tcilatex}
\begin{document}

\title{Non-Archimedean Probability }
\author{Vieri Benci\thanks{
Dipartimento di Matematica Applicata, Universit\`a degli Studi di Pisa, Via
F. Buonarroti 1/c, Pisa, ITALY and Department of Mathematics, College of
Science, King Saud University, Riyadh, 11451, SAUDI ARABIA. e-mail: \texttt{%
benci@dma.unipi.it}} \and Leon Horsten\thanks{
Department of Philosophy, University of Bristol, 43 Woodland Rd, BS81UU
Bristol, UNITED KINGDOM. e-mail: \texttt{leon.horsten@bristol.ac.uk}} \and %
Sylvia Wenmackers\thanks{
Faculty of Philosophy, University of Groningen, Oude Boteringestraat 52,
9712 GL Groningen, THE NETHERLANDS. e-mail: \texttt{s.wenmackers@rug.nl}} }
\maketitle

\begin{abstract}
We propose an alternative approach to probability theory closely related to
the framework of \textit{numerosity} theory: non-Archimedean probability
(NAP). In our approach, unlike in classical probability theory, all subsets
of an infinite sample space are measurable and zero- and unit-probability
events pose no particular epistemological problems. We use a non-Archimedean
field as the range of the probability function. As a result, the property of
countable additivity in Kolmogorov's axiomatization of probability is
replaced by a different type of infinite additivity.

\medskip


\medskip


\end{abstract}

\tableofcontents

%%%%%%%%%%%%%%%%%%%%%%%%%%%%%%

%%%%%%%%%%%%%%%%%%%%%%%%%%%%%%

\section{Introduction}

\label{sec:introduction} Kolmogorov's classical axiomatization embeds
probability theory into measure theory: it takes the domain and the range of
a probability function to be standard sets and employs the classical
concepts of limit and continuity. Kolmogorov starts from an elementary
theory of probability ``in which we have to deal with only a finite number
of events'' \cite[p.~1]{Kolmogorov 1933}. We will stay close to his axioms
for the case of finite sample spaces, but critically investigate his
approach in the second chapter of \cite{Kolmogorov 1933} dealing with the
case of ``an infinite number of random events''. There, Kolmogorov
introduces an additional axiom, the Axiom of Continuity. Together with the
axioms and theorems for the finite case (in particular, the addition
theorem, now called `finite additivity', FA), this leads to the generalized
addition theorem, called `$\sigma $-additivity' or `countable additivity'
(CA) in the case where the event space is a Borel field (or $\sigma $%
-algebra, in modern terminology). Some problems---such as a fair lottery on $%
\mathbb{N}$ or $\mathbb{Q}$---cannot be modeled within Kolmogorov's
framework. Weakening additivity to finite additivity, as de Finetti
advocated \cite{de Finetti 1974}, solves this but introduces strange
consequences of its own \cite{Kadane etal 1986}.

Nelson developed an alternative approach to probability theory based on
non-Archimedean, hyperfinite sets as the domain and range of the probability
function \cite{Nelson1987}. His framework has the benefit of making
probability theory on infinite sample spaces equally simple and
straightforward as the corresponding theory on finite sample spaces; the
appropriate additivity property is hyperfinite additivity. The drawback of
this elegant theory is that it does not apply \emph{directly} to
probabilistic problems on standard infinite sets, such as a fair lottery on $%
\mathbb{N}$, $\mathbb{Q}$, $\mathbb{R}$, or $2^\mathbb{N}$.

In the current paper, we develop a third approach to probability theory,
called non-Archimedean probability (NAP) theory. We formulate new desiderata
(axioms) for a concept of probability that is able to describe the case of a
fair lottery on $\mathbb{N}$, as well as other cases where infinite sample
spaces are involved. As such, the current article generalizes the solution
to the infinite lottery puzzle presented in \cite{Wenmackers and Horsten
2010}. Within NAP-theory, the domain of the probability function can be the
full powerset of any standard set from applied mathematics, whereas the
general range is a non-Archimedean field. We investigate the consistency of
the proposed axioms by giving a model for them. We show that our theory can
be understood in terms of a novel formalization of the limit and continuity
concept (called `$\Lambda$-limit' and `non-Archimedean continuity',
respectively). Kolmogorov's CA is replaced by a different type of infinite
additivity. In the last section, we give some examples. For fair lotteries,
the probability assigned to an event by NAP-theory is directly proportional
to the `numerosity' of the subset representing that event \cite{BDN2003a}.

There is a philosophical companion article \cite{BenciHorstenWenm} to this
article, in which we dissolve philosophical objections against using
infinitesimals to model probability on infinite sample spaces. We regard all
mentioned frameworks for probability theory as mathematically correct
theories (\textit{i.e.}~internally consistent), with a different scope of
applicability. Exploring the connections between the various theories helps
to understand each of them better. We agree with Nelson \cite{Nelson1987}
that the infinitesimal probability values should not be considered as an
intermediate step---a method to arrive at the answer---but rather as the
final answer to probabilistic problems on infinite sample spaces. Approached
as such, NAP-theory is a versatile tool with epistemological advantages over
the orthodox framework.

\subsection{Some notation}

Here we fix the notation used in this paper: we set

\begin{itemize}
\item $\mathbb{N}=\left\{ 1,2,3,.....\right\} $ is the set of positive
natural numbers;

\item $\mathbb{N}_{0}=\left\{ 0,1,2,3,.....\right\} $ is the set of natural
numbers;

\item if $A$ is a set, then $|A|$ will denote its cardinality;

\item if $A$ is a set, $\mathcal{P}(A)$ is the set of the parts of $A,$ $%
\mathcal{P}_{fin}(A)$ is the set of \textit{finite} parts of $A$;

\item for any set $A,$ $\chi _{A}$ will denote the characteristic function
of $A$, namely 
\begin{equation*}
\chi _{A}(\omega )=\left\{ 
\begin{array}{cc}
1 & if\ \ \omega \in A \\ 
&  \\ 
0 & if\ \ \omega \notin A ;%
\end{array}
\right.
\end{equation*}

\item if $F$ is an ordered field, and $a,b\in F,$ then we set 
\begin{eqnarray*}
\left[ a,b\right] _{F} &:&=\left\{ x\in F\ |\ a\leq x\leq b\right\} \\
\left[ a,b\right) _{F} &:&=\left\{ x\in F\ |\ a\leq x<b\right\} ;
\end{eqnarray*}

\item if $F$ is an ordered field, and $F\supseteq \mathbb{R}$, then $F$ is
called a superreal field;

\item for any set $\mathfrak{\mathcal{D}},$ $\mathfrak{F}\left( \mathfrak{%
\mathcal{D}},\mathbb{R}\right) $ will denote the (real) algebra of functions 
$u:\mathfrak{\mathcal{D}}\rightarrow \mathbb{R}$ equipped with the following
operations: for any $u, v \in \mathfrak{F}\left( \mathfrak{\mathcal{D}},%
\mathbb{R}\right)$, for any $r\in \mathbb{R}$, and for any $x \in \mathcal{D}
$: 
\begin{eqnarray*}
(u+v)(x) &=&u(x)+v(x), \\
\left( ru\right) (x) &=&ru(x), \\
(u\cdot v)(x) &=&u(x)\cdot v(x);
\end{eqnarray*}

\item if $F$ is a non-Archimedean field then we set 
\begin{equation*}
x\sim y \Leftrightarrow x-y\ \ \text{is infinitesimal;}
\end{equation*}
in this case we say that \ $x$ and $y$ are infinitely close;

\item if $F$ is a non-Archimedean superreal field and $\xi \in F$ is bounded
then $st(\xi )$ denotes the unique real number $x$ infinitely close to $\xi$.
\end{itemize}

\section{Kolmogorov's probability theory}

\subsection{Kolmogorov's axioms}

Classical probability theory is based on Kolmogorov's axioms (KA) \cite%
{Kolmogorov 1933}. We give an equivalent formulation of KA, using $P_{KA}$
to indicate a probability function that obeys these axioms. The sample space 
$\Omega$ is a set whose elements represent elementary events.

\begin{center}
{\Large Axioms of Kolmogorov}
\end{center}

\begin{itemize}
\item (K0) \textbf{Domain and range. }The events are the elements of $%
\mathfrak{A}$, a $\sigma$-algebra over $\Omega$,\footnote{$\mathfrak{A}$ is
a $\sigma$-algebra over $\Omega$ if and only if $\mathfrak{A} \subseteq 
\mathcal{P}\left( \Omega \right)$ such that $\mathfrak{A}$ is closed under
complementation, intersection, and countable unions. $\mathfrak{A}$ is
called the `event algebra' or `event space'.} and the probability is a
function 
\begin{equation*}
P_{KA}:\mathfrak{A} \rightarrow \mathbb{R}
\end{equation*}

\item (K1) \textbf{Positivity. }$\forall A\in \mathfrak{A},$ 
\begin{equation*}
P_{KA}(A)\geq 0
\end{equation*}

\item (K2) \textbf{Normalization.} 
\begin{equation*}
P_{KA}(\Omega )=1
\end{equation*}

\item (K3) \textbf{Additivity. }$\forall A, B \in \mathfrak{A}$ such that $%
A\cap B=\varnothing$, 
\begin{equation*}
P_{KA}(A\cup B)=P_{KA}(A)+P_{KA}(B)
\end{equation*}

\item (K4)\ \textbf{Continuity. }\textit{Let }%
\begin{equation*}
A=\dbigcup\limits_{n\in \mathbb{N}}A_{n}
\end{equation*}%
\textit{where }$\forall n\in \mathbb{N},A_{n}\subseteq A_{n+1}\subseteq 
\mathfrak{A}$\textit{; then} 
\begin{equation*}
P_{KA}(A)=\ \underset{n\in \mathbb{N}}{\sup }\ P_{KA}(A_{n})
\end{equation*}
\end{itemize}

We will refer to the triple $(\Omega ,\mathfrak{A},P_{KA})$ as a \textit{%
Kolmogorov Probability space}.

\begin{remark}
If the sample space is finite then it is sufficient to define a normalized
probability function on the elementary events, namely a function 
\begin{equation*}
p:\Omega \rightarrow \mathbb{R}
\end{equation*}
with 
\begin{equation*}
\sum_{\omega \in \Omega }p(\omega )=1
\end{equation*}
In this case the probability function 
\begin{equation*}
P_{KA}:\mathcal{P}\left( \Omega \right) \rightarrow \left[ 0,1\right] _{%
\mathbb{R}}
\end{equation*}
is defined by 
\begin{equation}
P_{KA}(A)=\sum_{\omega \in A}p(\omega )  \label{liu}
\end{equation}
and KA are trivially satisfied. Unfortunately, eq.~(\ref{liu}) cannot be
generalized to the infinite case. In fact, if the sample space is $\mathbb{R}
$, an infinite sum might yield a result in $[0,1]_\mathbb{R}$ only if $%
p(\omega)\neq 0$ for at most a denumerable number of $\omega \in A$.%
\footnote{%
Similarly, in classical analysis the sum of an uncountable sequence is
undefined.} In a sense, the Continuity Axiom (K4) replaces eq.~(\ref{liu})
for infinite sample spaces.
\end{remark}

\subsection{Problems with Kolmogorov's axioms\label{sec:PK}}

Kolmogorov uses $[0,1]_{\mathbb{R}}$ as the range of $P_{KA}$, which is a
subset of an ordered field and thus provides a good structure for adding and
multiplying probability values, as well as for comparing them. However, this
choice for the range of $P_{KA}$ in combination with the property of
Countable Additivity (which is a consequence of Continuity) may lead to
problems in cases with infinite sample spaces.

\bigskip

\noindent\textbf{Non-measurable sets in }$\mathcal{P}\left( \Omega \right)$
A peculiarity of KA is that, in general $\mathfrak{A}\neq \mathcal{P}\left(
\Omega \right)$. In fact, it is well known that there are (probability)
measures (such as the Lebesgue measure on $\left[ 0,1\right] $) which cannot
be defined for all the sets in $\mathcal{P}\left( \Omega\right)$. Thus,
there are sets in $\mathcal{P}\left( \Omega \right) $ which are not events
(namely elements of $\mathfrak{A}$) even when they are the union of
elementary events in $\mathfrak{A}$.\footnote{%
For example, if $\Omega =\left[ 0,1\right] _{\mathbb{R}}\ $and $P_{KA}$ is
given by the Lebesgue measure, then all the singletons $\left\{ x\right\}$
are measurable, but there are non-measurable sets; namely the union of
events might not be an event.}

\bigskip

\noindent\textbf{Interpretation of }$P_{KA}=0$\textbf{\ and }$P_{KA}=1$

In Kolmogorov's approach to probability theory, there are situations such
that:

\begin{equation}
P_{KA}(A_{j})=0,\ \ j\in J  \label{eq:zero}
\end{equation}%
and 
\begin{equation}
P_{KA}\left( \bigcup\limits_{j\in J}A_{j}\right) = 1.  \label{eq:one}
\end{equation}

\noindent This situation is very common when $J$ is not denumerable. It
looks as if eq.~(\ref{eq:zero}) states that each event $A_{j}$ is
impossible, whereas eq.~(\ref{eq:one}) states that one of them will occur
certainly. This situation requires further epistemological reflection.
Kolmogorov's probability theory works fine as a mathematical theory, but the 
\textit{direct} interpretation of its language leads to counterintuitive
results such as the one just described. An obvious solution is to interpret
probability 0 as `very unlikely' (rather than simply as `impossible'), and
to interpret probability 1 as `almost surely' (instead of `absolutely
certain'). Yet, there is a philosophical price to be paid to avoid these
contradictions: the correspondence between mathematical formulas and reality
is now quite vague---just how probable is `very likely' or `almost
surely'?---and far from intuition.

\bigskip

\noindent \textbf{Fair lottery on }$\mathbb{N}$

We may observe that the choice $[0,1]_{\mathbb{R}}$ as the range of the
probability function is neither necessary to describe a fair lottery in the
finite case, nor sufficient to describe one in the infinite case.

\begin{itemize}
\item For a fair finite lottery, the unit interval of $\mathbb{R}$ is not
necessary as the range of the probability function: the unit interval of $%
\mathbb{Q}$ (or, maybe some other denumerable subfield of $\mathbb{R}$)
suffices.

\item In the case of a fair lottery on $\mathbb{N}$, $[0,1]_{\mathbb{R}}$ is
not sufficient as the range: it violates our intuition that the probability
of any set of tickets can be obtained by adding the probabilities of all
individual tickets.
\end{itemize}

Let us now focus on the fair lottery on $\mathbb{N}$. In this case, the
sample space is $\Omega =\mathbb{N}$ and we expect the domain of the
probability function to contain all the singletons of $\mathbb{N}$,
otherwise there would be `tickets' (individual, natural numbers) whose
probability is undefined. Yet, we expect them to be defined and to be equal
in a fair lottery. Indeed, we expect to be able to assign a probability to
any possible combination of tickets. This assumption implies that the event
algebra should be $\mathcal{P}\left( \Omega \right)$. Moreover, we expect to
be able to calculate the probability of an arbitrary event by a process of
summing over the individual tickets.

This leads us to the following conclusions. First, if we want to have a
probability theory which describes a fair lottery on $\mathbb{N}$, assigns a
probability to all singletons of $\mathbb{N}$, and follows a generalized
additivity rule as well as the Normalization Axiom, the range of the
probability function has to be a subset of a non-Archimedean field. In other
words, the range has to include infinitesimals. Second, our intuitions
regarding infinite concepts are fed by our experience with their finite
counterparts. So, if we need to extrapolate the intuitions concerning finite
lotteries to infinite ones we need to introduce a sort of limit-operation
which transforms `extrapolations' into `limits'. Clearly, this operation
cannot be the limit of classical analysis. Since classical limits is
implicit in Kolmogorov's Continuity Axiom, this axiom must be revised in our
approach.

Motivated by the case study of a fair infinite lottery, at this point we
know which elements in Kolmogorov's classical axiomatization we do not
accept: the use of $[0,1]_{\mathbb{R}}$ as the range of the probability
function and the application of classical limits in the Continuity Axiom.
However, we have not specified an alternative to his approach: this is what
we present in the next sections.

\section{Non-Archimedean Probability}

\subsection{The axioms of Non-Archimedean Probability}

We will denote by $\mathfrak{F}\left( \mathfrak{\mathcal{P}}_{fin}(\Omega) , 
\mathbb{R} \right) $ the algebra of real functions defined on $\mathfrak{%
\mathcal{P}}_{fin}(\Omega)$.

\bigskip

\begin{center}
{\Large Axioms of Non-Archimedean Probability}
\end{center}

\begin{itemize}
\item (NAP0) \textbf{Domain and range. }\textit{The events are \textbf{all}
the elements of }$\mathcal{P}\left( \Omega \right) $\textit{\ and the
probability is a function} 
\begin{equation*}
P:\mathcal{P}\left( \Omega \right) \rightarrow \mathcal{R}
\end{equation*}
\textit{where }$\mathcal{R}$\textit{\ is a superreal field}

\item (NAP1) \textbf{Positivity.} $\forall A\in \mathcal{P}\left( \Omega
\right)$, 
\begin{equation*}
P(A)\geq 0
\end{equation*}

\item (NAP2) \textbf{Normalization.} $\forall A\in \mathcal{P}\left( \Omega
\right)$, 
\begin{equation*}
P(A)=1\Leftrightarrow A=\Omega
\end{equation*}

\item (NAP3) \textbf{Additivity.} $\forall A, B\in \mathcal{P}\left( \Omega
\right)$ such that $A\cap B=\varnothing$, 
\begin{equation*}
P(A\cup B)=P(A)+P(B)
\end{equation*}

\item (NAP4)\ \textbf{Non-Archimedean} \textbf{Continuity.} $\forall A, B\in 
\mathcal{P}\left( \Omega \right)$, with $B\neq \varnothing$, \textit{let }$%
P(A|B)$ \textit{\ denote the conditional probability, namely} 
\begin{equation}
P(A|B)=\frac{P(A\cap B)}{P(B)}.  \label{4+}
\end{equation}
Then

\begin{itemize}
\item $\forall \lambda \in \mathfrak{\mathcal{P}}_{fin}(\Omega ) \setminus
\varnothing,$ $P(A|\lambda )\in \mathbb{R}$\textit{\ };

\item \textit{there exists an algebra homomorphism} 
\begin{equation*}
J:\mathfrak{F}\left( \mathfrak{\mathcal{P}}_{fin}(\Omega ),\mathbb{R}\right)
\rightarrow \mathcal{R}
\end{equation*}
\textit{such that} $\forall A\in \mathfrak{\mathcal{P}}(\Omega )$ 
\begin{equation*}
P(A)=J\left( \varphi _{A}\right)
\end{equation*}
\textit{where} 
\begin{equation*}
\varphi _{A}(\lambda )=P(A|\lambda )\ \ \ \ for\ \ any\ \ \lambda \in 
\mathfrak{\mathcal{P}}_{fin}(\Omega).\footnote{%
In the remainder of this text, each occurence of $\lambda$ is to be
understood as referring to any $\lambda \in \mathfrak{\mathcal{P}}%
_{fin}(\Omega )$; $f(\lambda)$ will be used instead of $f(\cdot)$, where $f$
is a function on $\mathfrak{\mathcal{P}}_{fin}(\Omega)$.}
\end{equation*}
\end{itemize}
\end{itemize}

\bigskip

The triple $(\Omega ,P,J)$ will be called \textit{NAP-space} (or NAP-theory).

\bigskip

Now we will analyze the first three axioms and the fourth will be analyzed
in the next section.

The differences between (K0),...,(K3) and (NAP0),...,(NAP3) derive from
(NAP2). As consequence of this, we have that:

\bigskip

\begin{proposition}
\label{P}If (NAP0),...,(NAP3) holds, then

\begin{itemize}
\item (i) $\forall A\in \mathfrak{\mathcal{P}}(\Omega )$, $P(A)\in \left[ 0,1%
\right] _{\mathcal{R}}$

\item (ii) $\forall A\in \mathfrak{\mathcal{P}}(\Omega )$, $%
P(A)=0\Leftrightarrow A=\varnothing $

\item (iii) Moreover, there are sufficient conditions for $\mathcal{R}$ to
be non-Archimedean, such as:

\begin{itemize}
\item (a) $\Omega $ is countably infinite and the theory is fair, namely $%
\forall \omega ,\tau \in \Omega ,$ $P\left( \left\{ \omega \right\} \right)
=P\left(\left\{ \tau \right\} \right)$;

\item (b) $\Omega $ is uncountable.
\end{itemize}
\end{itemize}
\end{proposition}

\textbf{Proof.} Take $A\in \mathfrak{\mathcal{P}}(\Omega )$ and let $%
B=\Omega \setminus A$; then, by (NAP2) and (NAP3), 
\begin{equation*}
P(A)+P(B)=1;
\end{equation*}
then, since $P(B)\geq 0,$ $P(A)\leq 1$ and then (i) holds. Moreover, $%
P(A)=0\Leftrightarrow P(B)=1$ and hence, by (NAP2), $P(A)=0\Leftrightarrow
B=\Omega $ and so $P(A)=0\Leftrightarrow A=\varnothing .$ Now let us prove
(iii)(a) and assume that $\forall \omega \in \Omega ,$ $P\left(
\left\{\omega \right\} \right) =\varepsilon >0.$ Now we argue indirectly. If
the field $\mathcal{R}$ is Archimedean, then there exists $n\in \mathbb{N}$
such that $n\varepsilon >1;$ now let $A$ be a subset of $\Omega $ containing 
$n$ elements, then by (NAP3) it follows that $P\left( A\right) =nP\left(
\left\{\omega \right\} \right) =n\varepsilon >1$ and this fact contradicts
(NAP2); then $\mathcal{R}$ has to be non-Archimedean. Now let us prove
(iii)(b) and for every $n\in \mathbb{N}$ set 
\begin{equation*}
A_{n}=\left\{ \omega \in \Omega \ |\ 1/(n+1)<P\left( \left\{ \omega \right\}
\right) \leq 1/n\right\} .
\end{equation*}
By (NAP3) and (NAP2), it follows that each $A_{n}$ is finite, actually it
contains at most $n+1$ elements.

Now, again, we argue indirectly and assume that the field $\mathcal{R}$ is
Archimedean; in this case there are no infinitesimals and hence 
\begin{equation*}
\Omega =\dbigcup\limits_{n\in \mathbb{N}}A_{n}
\end{equation*}
and this contradicts the fact that we have assumed $\Omega $ to be
uncountable.

$\square $

\bigskip

\begin{remark}
\label{campo}In the axioms (NAP0),...,(NAP3), the field $\mathcal{R}$ is not
specified. This is not surprising since also in the Kolmogorov probability
the same may happen. For example, consider this simple example $\Omega
=\left\{ a,b\right\} ;$ $P_{KA}(\left\{ a\right\} )=1/\sqrt{2}$ and hence $%
P_{KA}(\left\{ b\right\} )=1-1/\sqrt{2}.$ In this case the natural field is $%
\mathbb{Q}(\sqrt{2}).$ However in Kolmogorov probability, since there is no
need to introduce infinitesimal probabilities, all these fields are
contained in $\mathbb{R}$ and hence it is simpler to take as range $\left[
0,1\right] _{\mathbb{R}}.$ We can make an analogous operation with NAP; this
will be done in section~\ref{errestar}.
\end{remark}

\subsection{Analysis of the fourth axiom\label{afa}}

If $A$ is a bounded subset of a non-Archimedean field then the supremum
might not exist; consider for example the set of all infinitesimal numbers.
Hence the axiom (K4) cannot hold in a non-Archimedean probability theory. In
this section, we will show an equivalent formulation of (K4) which can be
compared with (NAP4) and helps to understand the meaning of the latter.

\bigskip

\textbf{Conditional Probability Principle (CPP).} \textit{Let }$\Omega _{n}$%
\textit{\ be a family of events such that }$\Omega _{n}\subseteq
\Omega_{n+1} $\textit{\ and }$\Omega =\dbigcup\limits_{n\in \mathbb{N}%
}\Omega _{n};$\textit{\ then, eventually} 
\begin{equation*}
P_{KA}\left( \Omega _{n}\right) >0
\end{equation*}
\textit{and, for any event }$A,$\textit{\ we have that} 
\begin{equation*}
P_{KA}(A)=\ \underset{n\rightarrow \infty }{\lim }\ P_{KA}(A\ |\ \Omega
_{n}).
\end{equation*}

\bigskip

The following theorem shows that the Continuity Axiom (K4) is equivalent to
(CPP).

\begin{theorem}
Suppose that (K0),...,(K3) hold. (K4) holds if and only if the Conditional
Probability Principle holds.
\end{theorem}

\textbf{Proof:} Assume (K0),...,(K4) and let $\Omega _{n}$ be as in (CPP).
By (K2) and (K4) 
\begin{equation*}
\underset{n\in \mathbb{N}}{\sup }\ P_{KA}(\Omega _{n})=1
\end{equation*}
and so eventually $P_{KA}\left( \Omega _{n}\right) >0.$ Now take an event $A$%
. Since $P_{KA}(A\cap \Omega _{n})$ and $P_{KA}(\Omega _{n})$ are monotone
sequences, we have that 
\begin{eqnarray*}
\underset{n\rightarrow \infty }{\lim }\ P_{KA}(A\ |\ \Omega _{n}) &=&\ 
\underset{n\rightarrow \infty }{\lim }\ \frac{P_{KA}(A\cap \Omega _{n})}{%
P_{KA}(\Omega _{n})} \\
&=&\frac{\underset{n\in \mathbb{N}}{\sup }\ P_{KA}(A\cap \Omega _{n})}{%
\underset{n\in \mathbb{N}}{\sup }\ P_{KA}(\Omega_{n})}=\frac{P_{KA}(A)}{%
P_{KA}(\Omega )}=P_{KA}(A).
\end{eqnarray*}

Now assume (K0),....,(K3) and (CPP). Take any sequence $\Omega _{n}$ as in
(CPP); first we want to show that 
\begin{equation}
\underset{n\rightarrow \infty }{\lim }\ P_{KA}(\Omega _{n})=\ \underset{n\in 
\mathbb{N}}{\sup }\ P_{KA}(\Omega _{n})=1.  \label{rita}
\end{equation}
Take $\bar{n}$ such that $P_{KA}\left( \Omega _{\bar{n}}\right) >0;$ such a $%
\bar{n}$ exists since (CPP) holds. Then, using (CPP) again, we have 
\begin{equation*}
P_{KA}\left( \Omega _{\bar{n}}\right) =\ \underset{n\rightarrow \infty }{%
\lim }\ P_{KA}(\Omega _{\bar{n}}\ |\ \Omega _{n})=\frac{\underset{%
n\rightarrow \infty }{\lim }\ P_{KA}(\Omega _{\bar{n}}\cap \Omega _{n})}{%
\underset{n\rightarrow \infty }{\lim }\ P_{KA}(\Omega _{n})}=\frac{%
P_{KA}(\Omega _{\bar{n}})}{\underset{n\rightarrow \infty }{\lim }\
P_{KA}(\Omega _{n})}.
\end{equation*}
Since $P_{KA}\left( \Omega _{\bar{n}}\right) >0,\ $eq.~(\ref{rita}) follows.

Now let $A_{n}$ be a sequence as in (K4) and set 
\begin{equation*}
\Omega _{n}=\left( \Omega \setminus A\right) \cup A_{n}.
\end{equation*}
Then $\Omega _{n}$ and $A$ satisfies the assumptions of (CPP). So, by eq.~(%
\ref{rita}), we have that 
\begin{eqnarray*}
P_{KA}(A) &=&\ \underset{n \rightarrow \infty }{\lim }\ P_{KA}(A\ |\
\Omega_{n})=\frac{\underset{n \rightarrow \infty }{\lim }\ P_{KA}(A\cap
\Omega _{n})}{\underset{n \rightarrow \infty }{\lim }\ P_{KA}(\Omega _{n})}
\\
&=&\frac{\underset{n \rightarrow \infty }{\lim }\ P_{KA}(A_{n})}{1}=\ 
\underset{n \in \mathbb{N}}{\sup }\ P_{KA}(A_{n}).
\end{eqnarray*}

$\square $

\bigskip

So (CPP) is equivalent to (K4) and it has a form which can be compared with
(NAP4). Both (CPP) and (NAP4) imply that the knowledge of the conditional
probability relative to a suitable family of sets provides the knowledge of
the probability of the event. In the Kolmogorovian case, we have that 
\begin{equation}
P_{KA}(A)=\ \underset{n\rightarrow \infty }{\lim }\ P_{KA}(A\ |\ \Omega _{n})
\label{PAL}
\end{equation}
and in the NAP case, we have that 
\begin{equation}
P(A)=J\left( P(A|\cdot )\right) .  \label{JJJ}
\end{equation}
If we compare these two equations, we see that we may think of $J\ $as a
particular kind of limit; this fact justifies the name \textit{%
Non-Archimedean Continuity} given to (NAP4). This point will be developed in
section~\ref{JL}.

\bigskip

\subsection{First consequences of the axioms}

Define a function 
\begin{equation*}
p:\Omega \rightarrow \mathcal{R}
\end{equation*}
as follows: 
\begin{equation*}
p(\omega )=P\left( \left\{ \omega \right\} \right)
\end{equation*}
We choose arbitrarily a point$\ \omega _{0}\in \Omega ,\ $and we define the 
\textit{weight }function as follows:\textit{\ } 
\begin{equation*}
w(\omega )=\frac{p(\omega )}{p(\omega _{0})}.
\end{equation*}

\begin{proposition}
The function $w$ takes its values in $\mathbb{R}$ and for any finite $%
\lambda $, the following holds: 
\begin{equation}
P(A\ |\ \lambda )=\frac{\sum_{\omega \in A\cap \lambda }w\left( \omega
\right) }{\sum_{\omega \in \lambda }w\left( \omega \right) }.  \label{aa}
\end{equation}
\end{proposition}

\textbf{Proof.} Take $\omega \in \Omega $ arbitrarily and set $r=P(\left\{
\omega \right\} |\left\{ \omega ,\omega _{0}\right\} ).$ By (NAP4), $r\in 
\mathbb{R}$ and, by the definition of conditional probability (see eq.~(\ref%
{4+})), we have that 
\begin{equation*}
r=\frac{p(\omega )}{p(\omega )+p(\omega _{0})}=\frac{w\left( \omega \right) 
}{w\left( \omega \right) +1}<1
\end{equation*}%
and hence 
\begin{equation*}
w\left( \omega \right) =\frac{r}{1-r}\in \mathbb{R}.
\end{equation*}

Eq.~(\ref{aa}) is a trivial consequence of the additivity and the definition
of $w:$ 
\begin{equation*}
P(A\ |\ \lambda )=\frac{\sum_{\omega \in A\cap \lambda }p\left( \omega
\right) }{\sum_{\omega \in \lambda }p\left( \omega \right) }=\frac{%
\sum_{\omega \in A\cap \lambda }p(\omega _{0})w\left( \omega \right) }{%
\sum_{\omega \in \lambda }p(\omega _{0})w\left( \omega \right) }=\frac{%
\sum_{\omega \in A\cap \lambda }w\left( \omega \right) }{\sum_{\omega \in
\lambda }w\left( \omega \right) }.
\end{equation*}

$\square $

We recall that $\chi _{\lambda }$ denotes the characteristic function of $%
\lambda $.

\begin{lemma}
$\forall \omega \in \Omega ,$ we have%
\begin{equation}
J\left( \chi _{\lambda }\left( \omega \right) \right) =1  \label{chi}
\end{equation}%
and%
\begin{equation}
J\left( \sum_{\omega \in \lambda }w\left( \omega \right) \right) =\frac{1}{%
p(\omega _{0})}  \label{pp}
\end{equation}
\end{lemma}

\textbf{Proof.} We have that 
\begin{eqnarray*}
\chi _{\lambda }\left( \omega \right) \left[ 1-\chi _{\lambda }\left( \omega
\right) \right] &=&0 \\
\chi _{\lambda }\left( \omega \right) +\left[ 1-\chi _{\lambda }\left(
\omega \right) \right] &=&1
\end{eqnarray*}%
then, setting $\xi =J\left( 1-\chi _{\lambda }\left( \omega \right) \right)$%
, we have that 
\begin{eqnarray*}
J\left( \chi _{\lambda }\left( \omega \right) \right) \cdot \xi &=&0 \\
J\left( \chi _{\lambda }\left( \omega \right) \right) +\xi &=&1
\end{eqnarray*}
then $J\left( \chi _{\lambda }\left( \omega \right) \right) \ $is either 1
or 0. We will show that $J\left( \chi _{\lambda }\left( \omega _{0}\right)
\right) =1.$

By eq.~(\ref{aa}), since $w(\omega _{0})=1,$ we have that 
\begin{equation*}
p(\omega _{0})=J\left( P(\left\{ \omega _{0}\right\} |\ \lambda )\right)
=J\left( \frac{w(\omega _{0})\chi _{\lambda }\left( \omega _{0}\right) }{%
\sum_{\omega \in \lambda }w\left( \omega \right) }\right) =\frac{J\left(
\chi _{\lambda }\left( \omega _{0}\right) \right) }{J\left( \sum_{\omega \in
\lambda }w\left( \omega \right) \right) }
\end{equation*}

By Prop. \ref{P} (ii), $p(\omega _{0})>0$, and hence $J\left( \chi _{\lambda
}\left( \omega _{0}\right) \right) =1$ and $J\left( \sum_{\omega \in \lambda
}w\left( \omega \right) \right) =1/p(\omega _{0})$

$\square $

\subsection{Infinite sums}

The NAP axioms allow us to generalize the notion of sum to infinite sets in
such a way that eq.~(\ref{liu}) and eq.~(\ref{aa}) hold also for infinite
sets.

If $F\subset \Omega $ is a finite set and $x_{\omega }$ $\in \mathbb{R\ }$
for $\omega \in \Omega $, using eq. (\ref{chi}), we have that 
\begin{equation*}
\sum_{\omega \in F}x_{\omega }=\sum_{\omega \in F}\left[ x_{\omega }J\left(
\chi _{\lambda }\left( \omega \right) \right) \right] =J\left( \sum_{\omega
\in F}x_{\omega }\chi _{\lambda }\left( \omega \right) \right) =J\left(
\sum_{\omega \in F\cap \lambda }x_{\omega }\right)
\end{equation*}

The last term makes sense also when $F$ is infinite (since $F \cap \lambda$
is finite). Hence, it makes sense to give the following definition:

\bigskip

\begin{definition}
\label{infsum}For any set $A\in \mathfrak{\mathcal{P}}(\Omega )$ and any
function $u:A\rightarrow \mathbb{R}$, we set 
\begin{equation*}
\sum_{\omega \in A}u(\omega )=J\left( \sum_{\omega \in A\cap
\lambda}u(\omega )\right).
\end{equation*}
\end{definition}

Using this notation, by eq.~(\ref{pp}), it follows that 
\begin{equation*}
\sum_{\omega \in \Omega }w\left( \omega \right) =\frac{1}{p(\omega _{0})}
\end{equation*}
and hence we get that 
\begin{equation*}
P(A)=P(A\ |\ \Omega )=\frac{\sum_{\omega \in A}w\left( \omega \right) }{%
\sum_{\omega \in \Omega }w\left( \omega \right) }=\frac{1}{p(\omega _{0})}
\sum_{\omega \in A}w\left( \omega \right).
\end{equation*}

Moreover, for any set $A$ and $B$, we have that 
\begin{equation*}
P(A|B)=\frac{P(A\cap B)}{P(B)}=\frac{\frac{1}{p(\omega _{0})}\sum_{\omega
\in A\cap B}w\left( \omega \right) }{\frac{1}{p(\omega _{0})}\sum_{\omega
\in B}w\left( \omega \right) }=\frac{\sum_{\omega \in A\cap B}w\left( \omega
\right) }{\sum_{\omega \in B}w\left( \omega \right) }.
\end{equation*}
This equation extends eq.~(\ref{aa}) when $\lambda $ is infinite. Moreover,
taking $B=\Omega $, we get 
\begin{equation*}
P(A)=\frac{1}{p(\omega _{0})}\sum_{\omega \in A}w\left( \omega \right).
\end{equation*}
This equation extends eq.~(\ref{liu}) when $A$ is infinite since $p\left(
\omega \right) =w\left( \omega \right) p(\omega_{0}).$

The next proposition replaces $\sigma $-additivity:\footnote{%
Because it also holds for non-denumerably infinite sample spaces, this
proposition encapsulates what some philosophers have called `perfect
additivity'; see \textit{e.g.}~\cite[Vol.~1, p.~118]{de Finetti 1974}.}

\begin{proposition}
\textit{Let } 
\begin{equation*}
A=\dbigcup\limits_{j\in I}A_{j}
\end{equation*}
\textit{where }$I\ $is a family of indices of any cardinality and $A_{j}\cap
A_{k}=\varnothing $ for $j\neq k$\textit{; then} 
\begin{equation*}
P(A)=J(\sigma )
\end{equation*}
where 
\begin{equation*}
\sigma \left( \lambda \right) :=\sum_{j\in I}P\left( A_{j}|\lambda \right)
\end{equation*}
\end{proposition}

\textbf{Proof}. Since $\lambda $ is finite, $P\left( A_{j}|\lambda \right) $
can be computed just by making finite sums: 
\begin{equation*}
\sum_{j\in I}P\left( A_{j}|\lambda \right) =\frac{\sum_{j\in I}\sum_{\omega
\in A_{j}\cap \lambda }w\left( \omega \right) }{\sum_{\omega \in \lambda
}w\left( \omega \right) }=\frac{\sum_{\omega \in A\cap \lambda }w\left(
\omega \right) }{\sum_{\omega \in \lambda }w\left( \omega \right) }=P\left(
A|\lambda \right).
\end{equation*}
So we have 
\begin{equation*}
J(\sigma )=J\left( \sum_{j\in I}P\left( A_{j}|\lambda \right)
\right)=J\left( P\left( A|\lambda \right) \right) =P\left( A\right) .
\end{equation*}

$\square $

\section{NAP-spaces and $\Lambda $-limits}

In this section, we will show how to construct NAP-spaces. In particular
this construction shows that the NAP-axioms are not contradictory. Also, we
will introduce the notion of $\Lambda $-limit which will be useful in the
applications.

\subsection{Fine ideals}

Before constructing NAP-spaces, we give the definition and some properties
of fine ideals.

\begin{definition}
An ideal $I$ in the algebra $\mathfrak{F}\left( \mathfrak{\mathcal{P}}%
_{fin}(\Omega ),\mathbb{R}\right) $ is called fine\footnote{%
The name \textit{fine} ideal has been chosen since the maximal ultrafilter 
\begin{equation*}
\mathcal{U}:=\left\{ \varphi ^{-1}(0)\ |\ \varphi \in I\right\}
\end{equation*}
is a \textit{fine} ultrafilter.} if it is maximal and if for any $\omega \in
\Omega ,$ $1-\chi _{\lambda }\left( \omega \right) \in I.$
\end{definition}

\begin{proposition}
If $\Omega $ is an infinite set, then $\mathfrak{F}\left( \mathfrak{\mathcal{%
P}}_{fin}(\Omega ),\mathbb{R}\right) $ contains a fine ideal.
\end{proposition}

\textbf{Proof.} We set 
\begin{equation*}
I_{0}=\left\{ \varphi \in \mathfrak{F}\left( \mathfrak{\mathcal{P}}%
_{fin}(\Omega ),\mathbb{R}\right) \ |\ \exists \lambda _{0}\in \mathfrak{%
\mathcal{P}}_{fin}(\Omega )\in \Omega ,\forall \lambda \supseteq
\lambda_{0},\ \varphi \left( \lambda \right) =0\right\}.
\end{equation*}

It is easy to see that $I_{0}$ is an ideal; in fact:

- if $\forall \lambda \supseteq \lambda _{0},\ \varphi \left( \lambda
\right) =0\ $ and $\forall \lambda \supseteq \mu _{0},\ \psi \left( \lambda
\right) =0,$ then $\forall \lambda \supseteq \lambda _{0}\cup \mu _{0},\
\left( \varphi +\psi \right) \left( \lambda \right) =0$ and hence $\varphi
+\psi \in I_{0}$;

- if $\forall \lambda \supseteq \lambda _{0},\ \varphi \left( \lambda
\right) =0,$ then, $\forall \psi \in \mathfrak{F}\left( \mathfrak{\mathcal{P}%
}_{fin}(\Omega ),\mathbb{R}\right) ,$ we have that $\forall \lambda
\supseteq \lambda _{0},\ \varphi \left( \lambda \right) \cdot \psi \left(
\lambda \right) =0$ and hence $\varphi \cdot \psi \in I_{0}$.

Moreover, $1-\chi _{\lambda }\left( \omega \right) \in I_{0}$ since $%
1-\chi_{\lambda }\left( \omega \right) =0$ $\forall \lambda \supseteq
\lambda_{0}:=\left\{ \omega \right\} .$

The conclusion follows taking a maximal ideal $I$ containing $I_{0}$ which
exists by Krull's theorem.

$\square $

\begin{proposition}
Let $(\Omega ,P,J)$ be a NAP-space; then $\ker \left( J\right) $ is a fine
ideal.
\end{proposition}

\textbf{Proof.} Since $\mathcal{R}$ is a field, by elementary algebra it
follows that $\ker \left( J\right) $ is a prime ideal: but all prime ideals
in a real algebra of functions are maximal. Thus $\ker \left( J\right) $ is
a maximal ideal and by eq.~(\ref{chi}), it follows that it is fine.

$\square $

\subsection{Construction of NAP-spaces\label{CN}}

In the previous section, we have seen that, given a NAP-space $(\Omega ,P,J)$
(with $\Omega $ infinite), it is possible to define the weight function $%
w:\Omega \rightarrow \mathbb{R}^{+}$ and a fine ideal $I\subset \mathfrak{F}%
\left( \mathfrak{\mathcal{P}}_{fin}(\Omega ),\mathbb{R}\right) .$ In this
section, we will show that also the converse is possible, namely that in
order to define a NAP-space it is sufficient to assign

\begin{itemize}
\item the sample space $\Omega$;

\item a weight function $w:\Omega \rightarrow \mathbb{R}^{+}$;

\item a fine ideal $I$ in the algebra $\mathfrak{F}\left( \mathfrak{\mathcal{%
P}}_{fin}(\Omega ),\mathbb{R}\right) $.
\end{itemize}

The weight function allows to define the conditional probability of an event 
$A$ with respect to an event $\lambda \in \mathcal{P}_{fin}\left( \Omega
\right) $ according to the following formula 
\begin{equation*}
P(A\ |\ \lambda )=\frac{\sum_{\omega \in A\cap \lambda }w\left( \omega
\right) }{\sum_{\omega \in \lambda }w\left( \omega \right) }.
\end{equation*}

The fine ideal $I\ $allows to define an ordered field 
\begin{equation}
\mathcal{R}_{I}:=\frac{\mathfrak{F}\left( \mathfrak{\mathcal{P}}%
_{fin}(\Omega ),\mathbb{R}\right) }{I}  \label{RR}
\end{equation}
and an algebra homomorphism 
\begin{equation*}
J_{I}:\mathfrak{F}\left( \mathfrak{\mathcal{P}}_{fin}(\Omega ),\mathbb{R}
\right) \rightarrow \mathcal{R}_{I}
\end{equation*}%
given by the canonical projection, namely 
\begin{equation*}
J_{I}\left( \varphi \right) =\left[ \varphi \right] _{I}\ \ \text{where\ \ } %
\left[ \varphi \right] _{I}:=\varphi +I.
\end{equation*}

The map $J_{I}$ allows to define an infinite sum as in Def.~\ref{infsum} and
to define the probability function as follows:%
\begin{equation}
P_{I}(A)=\frac{\sum_{\omega \in A}w\left( \omega \right) }{\sum_{\omega \in
\Omega }w\left( \omega \right) }.  \label{PP}
\end{equation}

Thus we have obtained the following theorem:

\begin{theorem}
Given $(\Omega ,w,I),$ the triple $(\Omega ,P_{I},\mathcal{R}_{I})$ defined
by eq.~(\ref{PP}) and eq.~(\ref{RR}) is a NAP-space, namely it satisfies the
axioms (NAP0),...,(NAP4).
\end{theorem}

\begin{definition}
$(\Omega ,P_{I},\mathcal{R}_{I})$ will be called the NAP-space produced by $%
(\Omega ,w,I).$
\end{definition}

\subsection{The $\Lambda $-property\label{LS}}

In section~\ref{CN}, we have seen that in order to construct a NAP, it is
sufficient to assign a triple $(\Omega ,w,I).$ However, it is not possible
to define $I$ explicitly since its existence uses Zorn's lemma and no
explicit construction is possible. In any case it is possible to choose $I\ $
in such a way that the NAP-theory satisfies some other properties which we
would like to include in the model. Some of these properties will be
described in section~\ref{SA} which deals with the applications. In this
section we will give a general strategy to include these properties in the
theory.

\begin{definition}
A family of sets $\Lambda \subset \mathfrak{\mathcal{P}}_{fin}(\Omega )$ is
called a directed set if

\begin{itemize}
\item if $\lambda _{1},\lambda _{2}\in \Lambda ,$ then $\exists \mu \in
\Lambda $ such that $\lambda _{1}\cup \lambda _{2}\subset \mu $;

\item the union of all the elements of $\Lambda $ gives $\Omega .$
\end{itemize}
\end{definition}

The notion of directed set allows to enunciate the following property.

\begin{definition}
\label{bla}Given a directed set $\Lambda ,$ we say that a NAP-space $%
(\Omega,P,\mathcal{R})$ satisfies the $\Lambda $-property, if, given $A,B\in 
\mathfrak{\mathcal{P}}(\Omega )$ such that 
\begin{equation*}
\forall \lambda \in \Lambda ,\ P(A\cap \lambda )=P(B\cap \lambda )
\end{equation*}
then 
\begin{equation*}
P(A)=P(B).
\end{equation*}
\end{definition}

Given any directed set $\Lambda \subset \mathfrak{\mathcal{P}}%
_{fin}(\Omega), $ it is easy to construct a NAP-space which satisfies the $%
\Lambda$-property. Given $\Lambda $, we define the ideal 
\begin{equation*}
I_{0,\Lambda }:=\left\{ \varphi \in \mathfrak{F}\left( \mathfrak{\mathcal{P}}%
_{fin}(\Omega ),\mathbb{R}\right) \ |\ \forall \lambda \in \Lambda ,\
\varphi \left( \lambda \right) =0\right\} ;
\end{equation*}
by Krull's theorem there exists a maximal ideal $I_{\Lambda }\supset
I_{0,\Lambda }.$ It is easy to check that $I_{\Lambda }$ is a fine ideal.
Then, given a weight function $w,$ we can consider the NAP-space produced by 
$(\Omega ,w,I_{\Lambda })$ and we have that:

\begin{theorem}
\label{TLprop}The NAP-space produced by $(\Omega ,w,I_{\Lambda })$ satisfies
the $\Lambda $-property.
\end{theorem}

\textbf{Proof.} Given $A,B\in \mathfrak{\mathcal{P}}(\Omega )$ as in def.~%
\ref{bla}, set 
\begin{equation*}
\varphi \left( \lambda \right) =P(A\cap \lambda )-P(B\cap \lambda );
\end{equation*}
then $\forall \lambda \in \Lambda ,\ \varphi \left( \lambda \right) =0,\ $
and hence $\varphi \in I_{0,\Lambda }\subset I_{\Lambda }$ and so $J(\varphi
\left( \lambda \right) )=0.$ Then we have: 
\begin{eqnarray*}
P(A)-P(B) &=&J(P(A|\lambda ))-J(P(B|\lambda )) \\
&=&J\left( \frac{P(A\cap \lambda )-P(B\cap \lambda )}{P(\lambda )}\right) = 
\frac{J(\varphi \left( \lambda \right) )}{J(P(\lambda ))}=0.
\end{eqnarray*}

Then the $\Lambda $-property holds.

$\square $

\subsection{The $\Lambda $-limit \label{JL}}

If we compare eq.~(\ref{PAL}) and eq.~(\ref{JJJ}), it makes sense to think
of $J$ as a particular kind of limit and to write eq.~(\ref{JJJ}) as
follows: 
\begin{equation*}
P(A)=\ \underset{\lambda \in \mathfrak{\mathcal{P}}_{fin}(\Omega );\ \lambda
\uparrow \Omega }{\lim }P(A|\lambda ).
\end{equation*}

More in general, we can define the following limit: 
\begin{equation*}
J(\varphi )=\ \underset{\lambda \in \mathfrak{\mathcal{P}}_{fin}(\Omega );\
\lambda \uparrow \Omega }{\lim }\ \varphi (\lambda),
\end{equation*}
where $\varphi \in \mathfrak{F}\left( \mathfrak{\mathcal{P}}_{fin}(\Omega ),%
\mathbb{R}\right) .$ The above limit is determined by the choice of the
ideal $I_{\Lambda }\subset \mathfrak{F}\left( \mathfrak{\mathcal{P}}%
_{fin}(\Omega ),\mathbb{R}\right) ,$ and, by Th.~\ref{TLprop}, it depends on
the values that $\varphi $ assumes on $\Lambda ;$ we can assume that $%
\varphi \in \mathfrak{F}\left( \Lambda ,\mathbb{R}\right) .$ Then, many
properties of this limit depend on the choice of $\Lambda ;\ $so we will
call it $\Lambda $-limit, and we will use the following notation 
\begin{equation*}
J(\varphi )=\ \lim_{\lambda \in \Lambda }\ \varphi (\lambda )
\end{equation*}
which is simpler and carries more information. Notice that the $\Lambda $%
-limit, unlike the usual limit, exists for any function $\varphi \in 
\mathfrak{F}\left( \Lambda ,\mathbb{R}\right) $ and that it takes its values
in the non-Archimedean field $\mathcal{R}_{I_{\Lambda }}$.

The following theorem shows some other properties of the $\Lambda $-limit.
These properties, except (v), are shared by the usual limit. However, (v)
and the fact that the $\Lambda $-limit always exists, make this limit quite
different from the usual one.

\begin{theorem}
\label{Tlim}Let $\varphi ,\psi :\Lambda \rightarrow \mathbb{R}$; then

\begin{itemize}
\item (i) \textit{If }$\varphi _{r}(\lambda )=r$\textit{\ is} \textit{%
constant}, \textit{then} 
\begin{equation*}
\lim_{\lambda \in \Lambda }\varphi _{r}(\lambda )=r
\end{equation*}

\item (ii) 
\begin{equation*}
\lim_{\lambda \in \Lambda }\varphi (\lambda )+\lim_{\lambda \in \Lambda}\psi
(\lambda )=\lim_{\lambda \in \Lambda }\left( \varphi (\lambda )+\psi
(\lambda )\right)
\end{equation*}

\item (iii) 
\begin{equation*}
\lim_{\lambda \in \Lambda }\varphi (\lambda )\cdot \lim_{\lambda \in \Lambda
}\psi (\lambda )=\lim_{\lambda \in \Lambda }\left( \varphi (\lambda )\cdot
\psi (\lambda )\right)
\end{equation*}

\item (iv) \textit{If }$\varphi (\lambda )$\textit{\ and }$\psi (\lambda )$ 
\textit{are} \textit{eventually} \textit{equal}, \textit{namely} $\exists
\lambda _{0}\in \Lambda :\ \forall \lambda \supset \lambda _{0},\ \varphi
(\lambda )=\psi (\lambda ),$ \textit{then} 
\begin{equation*}
\lim_{\lambda \in \Lambda }\varphi (\lambda )=\lim_{\lambda \in \Lambda}\psi
(\lambda )
\end{equation*}

\item (v) \textit{If }$\varphi (\lambda )$\textit{\ and }$\psi (\lambda )$ 
\textit{are} \textit{eventually} \textit{different}, \textit{namely} $%
\exists \lambda _{0}\in \Lambda :\ \forall \lambda \supset \lambda _{0},\
\varphi (\lambda )\neq \psi (\lambda ),$ \textit{then} 
\begin{equation*}
\lim_{\lambda \in \Lambda }\varphi (\lambda )\neq \lim_{\lambda \in \Lambda
}\psi (\lambda )
\end{equation*}

\item (vi) \textit{If, for any }$\lambda ,$\textit{\ }$\varphi (\lambda )$ 
\textit{\ has finite range, namely }$\varphi (\lambda )\in
\left\{r_{1},....,r_{n}\right\} ,\ $then 
\begin{equation*}
\lim_{\lambda \in \Lambda }\varphi (\lambda )=r_{j}
\end{equation*}%
for some $j\in \left\{ 1,....,n\right\} $.
\end{itemize}
\end{theorem}

\textbf{Proof.} (i) We have that 
\begin{equation*}
\lim_{\lambda \in \Lambda }\varphi (\lambda )=J(r\cdot 1)=r\cdot J(1)=r\cdot
1=r.
\end{equation*}
(ii) and (iii) are immediate consequences of the fact that $J$ is an
homomorphism.

(iv) Suppose that $\forall \lambda \supset \lambda _{0},\ \varphi (\lambda
)=\psi (\lambda );$ we set $\lambda _{0}=\left\{ \omega
_{1},...,\omega_{n}\right\} $ and 
\begin{equation*}
\zeta \left( \lambda \right) =\chi _{\lambda }(\omega _{1})\cdot .....\cdot
\chi _{\lambda }(\omega _{n}).
\end{equation*}
If $\lambda _{0}\setminus \lambda \neq \varnothing ,$ then $\zeta
\left(\lambda \right) =0$ since some of the $\chi _{\lambda }(\omega _{j})$
vanish; if $\lambda _{0}\setminus \lambda =\varnothing ,\ $then $\varphi
(\lambda )-\psi (\lambda )=0$ by our assumptions; therefore, 
\begin{equation}
\zeta \left( \lambda \right) \cdot \left[ \varphi (\lambda )-\psi (\lambda ) %
\right] =0.  \label{della}
\end{equation}
Moreover, by eq.~(\ref{chi}), we have that 
\begin{equation*}
J(\zeta )=J(\chi _{\lambda }(\omega _{1}))\cdot .....\cdot J(\chi _{\lambda
}(\omega _{n}))=1
\end{equation*}
and so, by eq.~(\ref{della}), 
\begin{eqnarray*}
\lim_{\lambda \in \Lambda }\varphi (\lambda )-\lim_{\lambda \in \Lambda
}\psi (\lambda ) &=&J\left( \varphi -\psi \right) =J\left( \varphi -\psi
\right) \cdot J(\zeta ) \\
&=&J\left( \left[ \varphi -\psi \right] \cdot \zeta \right) =J(0)=0.
\end{eqnarray*}

(v) We set 
\begin{equation*}
\theta (\lambda )=\left\{ 
\begin{array}{cc}
1 & if\ \ \lambda _{0}\setminus \lambda \neq \varnothing \\ 
&  \\ 
\frac{1}{\varphi (\lambda )-\psi (\lambda )} & if\ \ \lambda \supset
\lambda_{0} ;%
\end{array}
\right.
\end{equation*}
then $\forall \lambda \supset \lambda _{0},$ 
\begin{equation*}
\left( \varphi (\lambda )-\psi (\lambda )\right) \cdot \theta (\lambda )=1
\end{equation*}
and hence, by (i) and (iv), 
\begin{eqnarray*}
1 &=&\lim_{\lambda \in \Lambda }\left[ \left( \varphi (\lambda )-\psi
(\lambda )\right) \cdot \theta (\lambda )\right] \\
&=&\lim_{\lambda \in \Lambda }\left( \varphi (\lambda )-\psi (\lambda
)\right) \cdot \lim_{\lambda \in \Lambda }\theta (\lambda ).
\end{eqnarray*}
From here, it follows that $\lim_{\lambda \in \Lambda }\left( \varphi
(\lambda )-\psi (\lambda )\right) \neq 0$ and by (ii) we get that $%
\lim_{\lambda \in \Lambda }\varphi (\lambda )\neq \lim_{\lambda \in \Lambda
}\psi (\lambda ).$

(vi) We have that 
\begin{equation*}
\left( \varphi (\lambda )-r_{1}\right) \cdot ......\cdot \left( \varphi
(\lambda )-r_{n}\right) =0;
\end{equation*}
then, taking the $\Lambda $-limit, 
\begin{equation*}
\left( \lim_{\lambda \in \Lambda }\varphi (\lambda )-r_{1}\right)
\cdot......\cdot \left( \lim_{\lambda \in \Lambda }\varphi
(\lambda)-r_{n}\right) =0;
\end{equation*}
and hence, there is a $j$ such that 
\begin{equation*}
\lim_{\lambda \in \Lambda }\varphi (\lambda )-r_{j}=0
\end{equation*}
and so $\lim_{\lambda \in \Lambda }\varphi (\lambda )=r_{j}.$

$\square $

If we use the notion of $\Lambda $-limit, def.~(\ref{infsum}) becomes more
meaningful; in this case in order to define an infinite sum $\sum_{\omega
\in A}u(\omega ),\ $we define the \textit{partial sum} as follows 
\begin{equation*}
\sum_{\omega \in A\cap \lambda }u(\omega ),\ \ \ \ \lambda \in \mathfrak{%
\mathcal{P}}_{fin}(\Omega )\ 
\end{equation*}
and then, we define the infinite sum as the $\Lambda $-limit of the partial
sums, namely 
\begin{equation*}
\sum_{\omega \in A}u(\omega )=\lim_{\lambda \in \Lambda }\sum_{\omega \in
A\cap \lambda }u(\omega ).
\end{equation*}

Moreover, the notion of $\Lambda $-limit provides also a meaningful
characterization of the field $\mathcal{R}_{I_{\Lambda }}:$ 
\begin{equation}
\mathcal{R}_{I_{\Lambda }}=\left\{ \lim_{\lambda \in \Lambda }\varphi
(\lambda )\ |\ \varphi \in \mathfrak{F}\left( \Lambda ,\mathbb{R}\right)
\right\}.  \label{erreL}
\end{equation}

Concluding, we have obtained the following `general strategy' for defining
NAP-spaces:

\bigskip

\textbf{General strategy. }In the applications, in order to define a
NAP-space we will assign

\begin{itemize}
\item the sample space $\Omega$;

\item a weight function $w:\Omega \rightarrow \mathbb{R}^{+}$;

\item a directed set $\Lambda \subseteq \mathfrak{\mathcal{P}}_{fin}(\Omega
) $ which provides a notion of $\Lambda $-limit and, via eq.~(\ref{erreL}),
the appropriate non-Archimedean field.
\end{itemize}

\bigskip

\subsection{The field $\mathbb{R}^{\ast }$\label{errestar}}

In this section, we will describe a non-Archimedean field $\mathbb{R}^{\ast} 
$ which contains the range of any non-Archimedean probability $P$ one may
wish to consider in applied mathematics.

To do this, we assume the existence of uncountable, non-accessible cardinal
numbers and, as usual, we will denote the smallest of them by $\kappa $. If
we assume the existence of $\kappa $, then there exists a nonstandard model $%
\mathbb{R}^{\ast }$ of cardinality $\kappa $ and $\kappa $-saturated. This
fact implies that it is unique up to isomorphisms. We refer to \cite[p.~194]%
{Keisler} for details.

Given a NAP-space $(\Omega ,P_{I},\mathcal{R}_{I})$, with $\Omega$ infinite,
we have that also $\mathcal{R}_{I}$ is a nonstandard model of $\mathbb{R}$
and if $|\Omega |\ <\kappa$, we have that $\mathcal{R}_{I}\subset \mathbb{R}%
^{\ast}$.

\bigskip

So using $\mathbb{R}^{\ast }$ and the notion of $\Lambda $-limit, axioms
(NAP0) and (NAP4) can be reformulated as follows:

\begin{itemize}
\item (NAP0)* \textbf{Domain and range. }\textit{The events are \textbf{all}
the elements of }$\mathcal{P}\left( \Omega \right) $\textit{\ and the
probability is a function} 
\begin{equation*}
P:\mathcal{P}\left( \Omega \right) \rightarrow \mathbb{R}^{\ast }
\end{equation*}%
\textit{where }$\mathbb{R}^{\ast }$ \textit{is the unique }$\kappa $\textit{%
-saturated} \textit{nonstandard model of }$\mathbb{R}$\textit{\ having
cardinality }$\kappa $\textit{.}

\item (NAP4)*\ \textbf{Non-Archimedean} \textbf{Continuity.} \textit{Let }$%
P(A|B),$ $B\neq \varnothing ,$\textit{\ denote the conditional probability,
then, }$P(A|\lambda )\in \mathbb{R}$ and 
\begin{equation*}
P(A)=\lim_{\lambda \in \Lambda }\ P(A|\lambda )
\end{equation*}%
\textit{for some directed set} $\Lambda \subset \mathfrak{F}\left( \mathfrak{%
\mathcal{P}}_{fin}(\Omega ),\mathbb{R}\right) $.
\end{itemize}

\begin{remark}
The previous remark shows some relation between NAP-theory and Nonstandard
Analysis. Actually, the relation is deeper than it appears here. In fact,
NAP could be constructed within a nonstandard universe based on the notion
of $\Lambda $-limit (see \cite{Bottazzi}). The idea to use Nonstandard
Analysis in probability theory is quite old and we refer to \cite%
{Keisler2002} and the references therein; these approaches differ from ours
since they use Nonstandard Analysis as a tool aimed at finding real-valued
probability functions. Another approach to probability related to
Nonstandard Analysis is due to Nelson \cite{Nelson1987}; however also his
approach is quite different from ours since it takes the domain of the
probability function to be a nonstandard set too.
\end{remark}

\bigskip

\section{Some applications\label{SA}}

\subsection{Fair lotteries and numerosities\label{num}}

\begin{definition}
If, $\forall \omega _{1},\omega _{2}\in \Omega ,$ $p(\omega_{1})=p(%
\omega_{2}),$ then the probability theory $(\Omega ,P,J)$ is called \emph{%
fair}.
\end{definition}

\bigskip

If $(\Omega ,P,J)$ is a fair lottery and $\Omega $ is finite, then, if we
set $\varepsilon _{0}=p(\omega _{0})=1/|\Omega |,\ $it turns out that, for
any set $A\in \mathfrak{\mathcal{P}}_{fin}(\Omega ),$ 
\begin{equation*}
|A|\ =\frac{P\left( A\right) }{\varepsilon _{0}}.
\end{equation*}

This remark suggests the following definition:

\bigskip

\begin{definition}
If $(\Omega ,P,J)$ is a fair lottery, the numerosity of a set $A \in 
\mathfrak{\mathcal{P}}(\Omega )$ is defined as follows: 
\begin{equation*}
\mathfrak{n}\left( A\right) =\frac{P\left( A\right) }{\varepsilon _{0}}
\end{equation*}
where $\varepsilon _{0}$ is the probability of an elementary event in $%
\Omega. $
\end{definition}

In particular if $A$ is finite, we have that $\mathfrak{n}\left(
A\right)=|A|.$ So the numerosity is the generalization to infinite sets of
the notion of ``number of elements of a set'' different from the Cantor
theory of infinite sets.

The theory of numerosity has been introduced in \cite{benci1995} and
developed in various directions in \cite{GR1996}, \cite{BDN2003}, and \cite%
{BDNF2006}. The definition above is an alternative way to introduce a
numerosity theory.

\bigskip

We now set 
\begin{equation*}
\mathbb{Q}^{\ast }=\left\{ \underset{\lambda \in \Lambda }{\lim }\ \varphi
(\lambda )\ |\ \forall \lambda \in \Lambda ,\ \varphi (\lambda )\in \mathbb{Q%
}\ \text{for some }\Lambda \mathcal{\ }\text{with }\left\vert \Lambda
\right\vert <\kappa \right\} .
\end{equation*}%
\bigskip

We will refer to $\mathbb{Q}^{\ast }$ as the field of hyperrational numbers.

\begin{proposition}
\label{FL}If $(\Omega ,P,J)$ is a fair lottery, and $A,B\subseteq \Omega $,
with $B\neq \varnothing$, then 
\begin{equation*}
P(A|B)\in \mathbb{Q}^{\ast }
\end{equation*}
\end{proposition}

\textbf{Proof}. We have that 
\begin{equation*}
P(A|B)=\frac{P(A\cap B)}{P(B)}=\underset{\lambda \in \mathfrak{\mathcal{P}}%
_{fin}(\Omega )}{\lim }\frac{P(A\cap B\cap \lambda )}{P(B\cap \lambda )}=%
\underset{\lambda \in \mathfrak{\mathcal{P}}_{fin}(\Omega )}{\lim }\frac{%
\left\vert A\cap B\cap \lambda \right\vert }{\left\vert B\cap \lambda
\right\vert }.
\end{equation*}
The conclusion follows from the fact that 
\begin{equation*}
\frac{\left\vert A\cap B\cap \lambda \right\vert }{\left\vert B\cap \lambda
\right\vert }\in \mathbb{Q}.
\end{equation*}

$\square $

\subsection{A fair lottery on $\mathbb{N}$}

Let us consider a fair lottery in which exactly one winner is randomly
selected from a countably infinite set of tickets. We assume that these
tickets are labeled by the (positive) natural numbers. Now let us construct
a NAP-space for such a lottery using the general strategy developed in
section~\ref{JL}.\footnote{%
This example has been discussed by multiple philosophers of probability,
including de Finetti \cite{de Finetti 1974}; at the end of this section, we
indicate how our solution relates to his approach. The solution presented
here rephrases the one given in \cite{Wenmackers and Horsten 2010} within
the more general NAP framework.}

We take 
\begin{equation*}
\Omega =\mathbb{N},
\end{equation*}
\begin{equation*}
w:\mathbb{N}\rightarrow \mathbb{R}^{+}\ \text{identically equal to\ one}
\end{equation*}%
and%
\begin{equation}
\Lambda _{\left[ n\right] }=\left\{ \lambda _{n}\in \mathfrak{\mathcal{P}}%
_{fin}(\Omega )\ |\ n\in \mathbb{N}\right\}  \label{laida}
\end{equation}
where 
\begin{equation*}
\lambda _{n}=\left\{ 1,2,3,.....,n\right\}.
\end{equation*}

In this case we have that, for every $A\in \mathfrak{\mathcal{P}}(\Omega )$ 
\begin{equation*}
P(A|\lambda _{n})=\frac{\left\vert A\cap \left\{ 1,....,n\right\}
\right\vert }{n}
\end{equation*}
and hence 
\begin{equation*}
P(A)=\lim_{\lambda _{n}\in \Lambda }\ \ \frac{\left\vert A\cap
\left\{1,....,n\right\} \right\vert }{n}.
\end{equation*}
Using the notion of numerosity as defined in section~\ref{num}, we set 
\begin{equation*}
\alpha :=\mathfrak{n(}\mathbb{N)=}\frac{1}{p(1)};
\end{equation*}
then the probability of any event $A\in \mathfrak{\mathcal{P}}(\Omega )$ can
be written as follows: 
\begin{equation*}
P(A)=\frac{\mathfrak{n(}A\mathbb{)}}{\alpha }
\end{equation*}
namely the probability of $A$ is ratio of the ``number'' of the elementary
events in $A$ and the total number of elements $\alpha .$

One of the properties which our intuition wants to be satisfied by the de
Finetti lottery is the `Asymptotic Limit Property', which relates the
non-Archimedean probability function $P$ to a classical limit (in as far as
the latter exists):

\begin{definition}
\textbf{Asymptotic Limit Property: }Let $A\in \mathcal{P}(\mathbb{N})$ be a
set which has an asymptotic limit, namely there exists $L\in \left[ 0,1%
\right] $ such that 
\begin{equation*}
\lim_{n\rightarrow \infty }\ \frac{\left\vert A\cap \left\{ 1,....,n\right\}
\right\vert }{n}=L.
\end{equation*}
We say that the Asymptotic Limit Property holds if we have 
\begin{equation}
P(A)\sim L.  \label{alp}
\end{equation}
\end{definition}

It is easy to check that

\begin{proposition}
The de Finetti lottery produced by $\left( \mathbb{N},1,I_{\Lambda _{\left[ n%
\right] }}\right) $ satisfies the Asymptotic Limit Property.
\end{proposition}

Now, we will consider additional properties which would be nice to have and
we will show how the choice of $\Lambda $ works. For example, the
probability of extracting an even number seems to be equal of extracting an
odd number; thus we must have 
\begin{equation}
P(\mathbb{E})=P(\mathbb{O})  \label{e1}
\end{equation}%
and since by (NAP2) and (NAP3), we have 
\begin{equation*}
P(\mathbb{E})+P(\mathbb{O})=1,
\end{equation*}%
it follows that 
\begin{equation}
P(\mathbb{E})=P(\mathbb{O})=\frac{1}{2}.  \label{e3}
\end{equation}

Now, let us compute for example $P(\mathbb{E}).$ We have that 
\begin{equation*}
P(\mathbb{E})=\frac{\mathfrak{n}(\mathbb{E})}{\alpha }
\end{equation*}
so we have to compute $\mathfrak{n}(\mathbb{E}):$ 
\begin{eqnarray*}
\mathfrak{n}(\mathbb{E}) &=&\lim_{\lambda _{n}\in \Lambda }\ \left\vert 
\mathbb{E}\cap \left\{ 1,....,n\right\} \right\vert \\
&=&\lim_{\lambda _{n}\in \Lambda }\ \left\vert \left\{ 2,4,6,....,2\cdot 
\left[ \frac{n}{2}\right] \right\} \right\vert \\
&=&\lim_{\lambda _{n}\in \Lambda }\ \left[ \frac{n}{2}\right]
=\lim_{\lambda_{n}\in \Lambda }\ \frac{n}{2}-\lim_{\lambda _{n}\in \Lambda
}\ c_{n}
\end{eqnarray*}
where $\left[ r\right] $ denotes the integer part of $r$ and 
\begin{equation*}
c_{n}=\left\{ 
\begin{array}{cc}
1/2 & if\ n\ is\ odd \\ 
0 & if\ n\ is\ even.%
\end{array}
\right.
\end{equation*}
Then, by Theorem~\ref{Tlim}(vi), $\lim_{\lambda _{n}\in \Lambda }c_{n}$ is
either $0$ or $1/2$, but this fact cannot be determined since we do not know
the ideal $I_{\Lambda }.$ However, if we think that this fact is relevant
for our model, we can follow the strategy suggested in section~\ref{LS} and
make a better choice of $\Lambda .$ If we choose a smaller $\Lambda ,$ we
have more information.

For example, we can replace the choice of eq.~(\ref{laida}) with the
following one: 
\begin{equation*}
\Lambda _{\left[ 2m\right] }:=\left\{ \left\{ 1,2,....,2m\right\} \ |\ m\in 
\mathbb{N}\right\}.
\end{equation*}

In this case we have that 
\begin{equation*}
\mathfrak{n}(\mathbb{E})=\lim_{\lambda _{n}\in \Lambda _{\left[ 2m\right]
}}\ \left[ \frac{n}{2}\right] =\lim_{\lambda _{n}\in \Lambda _{\left[ 2m%
\right] }}\ \frac{n}{2}=\frac{\alpha }{2}.
\end{equation*}
On the other hand, we can choose 
\begin{equation*}
\Lambda =\Lambda _{_{\left[ 2m-1\right] }}:=\left\{
\left\{1,2,....,2m-1\right\} \ |\ m\in \mathbb{N}\right\}
\end{equation*}
and in this case 
\begin{equation*}
\mathfrak{n}(\mathbb{E})=\lim_{\lambda \in \Lambda _{_{\left[ 2m-1\right]
}}}\ \left( \left[ \frac{n}{2}\right] -\frac{1}{2}\right) =\lim_{\lambda \in
\Lambda _{_{\left[ 2m-1\right] }}}\ \left( \frac{n-1}{2}\right) =\frac{%
\alpha -1}{2}.
\end{equation*}
Thus $P(\mathbb{E})=\frac{1}{2}$ or $\frac{1}{2}-\frac{1}{2\alpha }$
depending on the choice of $\Lambda .\ $Also, it is possible to prove that
any choice of $\Lambda \subseteq \Lambda _{\left[ n\right] }$ gives one of
these two possibilities.

\begin{remark}
The equality 
\begin{equation}
P(A)=L  \label{su}
\end{equation}%
cannot replace eq.~(\ref{alp}) for all the sets which have an asymptotic
limit; in fact take two sets $A$ and $B=A\cup F$ where $F$ is a finite set
(with $A\cap F=\varnothing $). Then, if $L$ is the asymptotic limit of $A,$
then it is also the asymptotic limit of $B$ since 
\begin{eqnarray*}
&&\lim_{n\rightarrow \infty }\ \frac{\left\vert B\cap
\left\{1,....,n\right\} \right\vert }{n} \\
&=&\lim_{n\rightarrow \infty }\ \frac{\left\vert A\cap
\left\{1,....,n\right\} \right\vert }{n}+\frac{\left\vert F\cap
\left\{1,....,n\right\} \right\vert }{n} \\
&=&L+0=L.
\end{eqnarray*}
On the other hand, by (NAP3) 
\begin{equation*}
P(B)=P(A)+P(F)
\end{equation*}
and by (NAP1), $P(F)>0\ $and hence $P(B)\neq P(A).$ Thus, it is not possible
that $P(A)=L\ \ $\textbf{and} $P(B)=L$.
\end{remark}

\bigskip

Even if eq.~(\ref{su}) cannot hold for all the sets, our intuition suggests
that in some cases it should be true and it would be nice if eq.~(\ref{su})
holds for a \textit{distinguished} family of sets. For example, if we have
eq.~(\ref{e1}) then eq.~(\ref{e3}) holds, and hence $P(\mathbb{E})$ and $P(%
\mathbb{O})$ have the probability equal to their asymptotic limit. So, the
following question arises naturally: is it possible to have a `de Finetti
probability space' produced by $\left\{ \mathbb{N},1,I_{\Lambda }\right\} $
such that 
\begin{equation*}
P(\mathbb{N}_{k})=\frac{1}{k}
\end{equation*}%
where 
\begin{equation*}
\mathbb{N}_{k}=\left\{ k,\ 2k,\ 3k,....,nk,....\right\} .
\end{equation*}%
The answer is yes; it is sufficient to choose 
\begin{equation}
\Lambda =\Lambda _{\left[ m!\right] }:=\left\{ \left\{ 1,....,m!\right\} \
|\ m\in \mathbb{N}\right\} .  \label{Lfatt}
\end{equation}%
In fact, in this case, we have 
\begin{eqnarray*}
P(\mathbb{N}_{k}) &=&\frac{\lim_{\lambda _{n}\in \Lambda _{\left[ m!\right]
}}\ \left\vert \mathbb{N}_{k}\mathbb{\cap }\left\{ 1,....,m!\right\}
\right\vert }{\alpha } \\
&=&\frac{\lim_{\lambda _{n}\in \Lambda _{\left[ m!\right] }}\ \left[ \frac{n%
}{k}\right] }{\alpha }=\frac{\lim_{\lambda _{n}\in \Lambda _{\left[ m!\right]
}}\ \frac{n}{k}}{\alpha }=\frac{1}{k}.
\end{eqnarray*}

More in general, if $\Lambda =\Lambda _{\left[ m!\right] },$ we can prove
that the sets 
\begin{equation*}
\mathbb{N}_{k,l}=\left\{ k-l,\ 2k-l,\ 3k-l,....,nk-l,....\right\} ,\
l=0,..,k-1
\end{equation*}
have probability $1/k,$ namely the same probability as their asymptotic
limit. This generalizes the situation which we have analyzed before where $%
\mathbb{E}=\mathbb{N}_{2}$ and $\mathbb{O}=\mathbb{N}_{2,1}.$

\begin{remark}
Our construction of a non-Archimedean probability $P,$ allows to construct
the following \emph{Archimedean }probability 
\begin{equation*}
P_{Arch}(A)=st(P(A))
\end{equation*}
where $st\left( \xi \right) $ denotes the standard part of $\xi ,$ namely
the unique standard number infinitely close to $\xi$. $P_{Arch}$ is defined
on \textbf{all} the subsets of $A$ and it is finitely additive and it
coincides with the asymptotic limit when it exists. Although we prefer a
theory based on non-Archimedean probabilities, we regard de Finetti's
reaction to the infinite lottery puzzle as an equally valid approach.%
\footnote{%
De Finetti was aware that probability could be treated as a non-Archimedean
quantity, but rejected this approach as ``a useless complication of
language'', which ``leads one to puzzle over `les infiniment petits'\ '' 
\cite[Vol.~2, p.~347]{de Finetti 1974}. He proposed to stay within an
Archimedean range, but to relax Kolmogorov's countable additivity to finite
additivity.} The construction of $P_{Arch}$ shows how the two approaches are
connected.
\end{remark}

\subsection{A fair lottery on $\mathbb{Q}$\label{FLQ}}

A fair lottery on $\mathbb{Q}$, by definition, is a NAP-space produced by $%
\left( \mathbb{Q},1,I\right) \ $for any arbitrary $I$; however, as in the
case of de Finetti lottery, we are allowed to require some additional
properties which appear natural to our intuition and then, we can inquire if
they are consistent.

For example, if we have two intervals $\left[ a_{0},b_{0}\right) _{\mathbb{Q}%
}\subset \left[ a_{1},b_{1}\right) _{\mathbb{Q}},$ we expect the conditional
probability to satisfy the following formula: 
\begin{equation}
P(\left[ a_{0},b_{0}\right) _{\mathbb{Q}}|\left[ a_{1},b_{1}\right) _{%
\mathbb{Q}})=\frac{b_{0}-a_{0}}{b_{1}-a_{1}}.  \label{CPLd}
\end{equation}

In fair lotteries, the probability is strictly related to the notion of
numerosity and the above formula is equivalent to the following one 
\begin{equation}
\mathfrak{n}(\left[ a,b\right) _{\mathbb{Q}})=\left( b-a\right) \cdot 
\mathfrak{n}(\left[ 0,1\right) _{\mathbb{Q}})  \label{linda}
\end{equation}
namely that the ``\textit{number of elements contained in an interval}'' is
proportional to its length.

Clearly, eq.~(\ref{CPLd}) follows from eq.~(\ref{linda}). In fact, 
\begin{equation*}
P(\left[ a,b\right) _{\mathbb{Q}})=\frac{\mathfrak{n}(\left[ a,b\right) _{%
\mathbb{Q}})}{\mathfrak{n}(\mathbb{Q})}=\left( b-a\right) \cdot \frac{%
\mathfrak{n}(\left[ 0,1\right) _{\mathbb{Q}})}{\mathfrak{n}(\mathbb{Q})}.
\end{equation*}
Then 
\begin{equation*}
P(\left[ a_{0},b_{0}\right) _{\mathbb{Q}}|\left[ a_{1},b_{1}\right) _{%
\mathbb{Q}})=\frac{P(\left[ a_{0},b_{0}\right) _{\mathbb{Q}})}{P(\left[
a_{1},b_{1}\right) _{\mathbb{Q}})}=\frac{b_{0}-a_{0}}{b_{1}-a_{1}}.
\end{equation*}
Also, it is easy to check that eq.~(\ref{linda}) follows from eq.~(\ref{CPLd}%
).

In order to prove that the property of eq.~(\ref{linda}) is consistent with
a NAP-theory, it is sufficient to find an appropriate family $\Lambda
\subset \mathcal{P}_{fin}(\mathbb{Q}).$

We will consider the family 
\begin{equation}
\Lambda _{\mathbb{Q}}:=\left\{ \mu _{n}\ |\ n=m!,\ m\in \mathbb{N}\right\}
\label{lambdaQ}
\end{equation}
with 
\begin{eqnarray}
\mu _{n} &=&\left\{ \frac{p}{n}\ |\ p\in \mathbb{Z},\ \left\vert \frac{p}{n}
\right\vert \leq n\right\}  \label{ln} \\
&=&\left\{ -n,-\frac{n^{2}-1}{n},....,-\frac{1}{n},0,\frac{1}{n},.....,\frac{%
n^{2}-1}{n},\ n\right\}.  \notag
\end{eqnarray}
In this case we have that:

\begin{proposition}
If $\Lambda =\Lambda _{\mathbb{Q}}$ is defined by eq.~(\ref{lambdaQ}), then
eq.~(\ref{linda}) holds.
\end{proposition}

\textbf{Proof}. We write\textbf{\ }$a$ and $b$ as fractions with the same
denominator: 
\begin{equation*}
a=\frac{p_{a}}{q};\ b=\frac{p_{b}}{q}.
\end{equation*}

Then, if you take $n$ sufficiently large (e.g. $n=m!,\ m\geq \max (a,b,q)$), 
$q$ divides $n$ and we have that 
\begin{equation*}
\left\vert \left[ a,b\right) _{\mathbb{Q}}\cap \mu _{n}\right\vert
=\left\vert \left\{ a,\ a+\frac{1}{n},\ a+\frac{2}{n},.....,b-\frac{1}{n}
\right\} \right\vert =\left( b-a\right) \cdot n.
\end{equation*}
From here, eq.~(\ref{linda}) easily follows.

$\square $

\bigskip

Let us compare the NAP-spaces produced by $\left( \mathbb{N},1,I_{\Lambda_{%
\mathbb{N}}}\right) \ $and $\left( \mathbb{Q},1,I_{\Lambda_{\mathbb{Q}%
}}\right)$, where we have set 
\begin{equation}
\Lambda _{\mathbb{N}}=\left\{ \mu _{n}\cap \mathbb{N}\ |\ \mu _{n}\in
\Lambda _{\mathbb{Q}}\right\} .  \label{LEnne}
\end{equation}

We want to show that this choice of $\Lambda _{\mathbb{N}}$ and $\Lambda _{%
\mathbb{Q}}\ $makes this two NAP-theories consistent in the sense described
below. First of all, notice that $\Lambda _{\mathbb{N}}$ defined by eq.~(\ref%
{LEnne}) coincides with $\Lambda _{\left[ m!\right] }$ defined by eq.~(\ref%
{Lfatt}) and that $\lambda _{n}=\mu _{n}\cap \mathbb{N}$.

If we denote by $P_{\mathbb{N}}$ and $P_{\mathbb{Q}}$ the respective
probabilities and we take $A\subset \mathbb{N}\subset \mathbb{Q},$ we have
that 
\begin{equation}
P_{\mathbb{N}}(A)=P_{\mathbb{Q}}(A\ |\ \mathbb{N}).  \label{PNA}
\end{equation}
In fact 
\begin{equation*}
P_{\mathbb{Q}}(A\ |\ \mathbb{N})=\frac{\mathfrak{n}\left( A\cap \mathbb{N}%
\right) }{\mathfrak{n}\left( \mathbb{N}\right) }=\frac{\mathfrak{n}\left(
A\right) }{\mathfrak{n}\left( \mathbb{N}\right) }=P_{\mathbb{N}}(A)
\end{equation*}
Moreover, if we set, as usual, $\alpha :=\mathfrak{n}\left( \mathbb{N}
\right) $, it is easy to check that 
\begin{eqnarray*}
\mathfrak{n}\left( \mathbb{Q}^{+}\right) &=&\alpha ^{2} \\
\mathfrak{n}\left( \mathbb{Q}\right) &=&2\alpha ^{2}+1.
\end{eqnarray*}
Thus, in this model, in a `lottery with rational numbers' the probability of
extracting a positive natural number is 
\begin{equation*}
P_{\mathbb{Q}}(\mathbb{N})=\frac{\mathfrak{n}\left( \mathbb{N}\right) }{%
\mathfrak{n}\left( \mathbb{Q}\right) }=\frac{\alpha }{2\alpha ^{2}+1}=\frac{%
1-\varepsilon }{2\alpha }
\end{equation*}
where $\varepsilon \ $is a positive infinitesimal.

\subsection{A fair lottery on $\mathbb{R}$}

Now let us consider a fair lottery on $\mathbb{R},$ namely a NAP-space
produced by $\left( \mathbb{R},1,I_{\Lambda _{\mathbb{R}}}\right) $ and let
us examine the other properties which we would like to have.\footnote{%
The problem of a fair lottery on a non-denumerable sample space is usually
presented as a fair lottery (or random darts throw) on the unit interval of $%
\mathbb{R}$ (or a darts board whose perimeter is indexed by this interval),
see \textit{e.g.}~\cite{Hajek2003}. The related problem of a random darts
throw on the unit square of $\mathbb{R}^{2}$ is considered in~\cite%
{BarthaHitchcock1999}.} Considering the example of the previous section, we
would like to have the analog of eq.~(\ref{CPLd}) just replacing $\mathbb{Q}$
with $\mathbb{R}$. This is not possible. Let us see why not. If in eq.~(\ref%
{CPLd}) we take $\left[ a_{0},b_{0}\right) _{\mathbb{R}}=\left[ 0,1\right) _{%
\mathbb{R}}\ $and $\left[ a_{1},b_{1}\right) _{\mathbb{R}}=\left[ 0,\sqrt{2}%
\right) _{\mathbb{R}},$ we would get 
\begin{equation*}
P\left( [0,1)_{\mathbb{R}}|[0,\sqrt{2})_{\mathbb{R}}\right) =\frac{1}{\sqrt{2%
}}
\end{equation*}%
and this fact is not possible by Prop.~\ref{FL} (in fact $\frac{1}{\sqrt{2}}%
\notin \mathbb{Q}^{\ast }$ by Th.~\ref{Tlim} (v) and the definition of $%
\mathbb{Q}^{\ast }$). However, we can require a weaker statement, namely
that, given two intervals $\left[ a_{0},b_{0}\right) _{\mathbb{R}}\subset %
\left[ a_{1},b_{1}\right) _{\mathbb{R}},$ 
\begin{equation}
P(\left[ a_{0},b_{0}\right) _{\mathbb{R}}|\left[ a_{1},b_{1}\right) _{%
\mathbb{R}})\sim \frac{b_{0}-a_{0}}{b_{1}-a_{1}}.  \label{tre}
\end{equation}%
Actually, in terms of numerosities, we can require that 
\begin{equation}
\forall a,b\in \mathbb{Q},\mathbb{\ }\mathfrak{n}(\left[ a,b\right) _{%
\mathbb{R}})=\left( b-a\right) \cdot \mathfrak{n}(\left[ 0,1\right) _{%
\mathbb{R}}  \label{uno}
\end{equation}%
\begin{equation}
\forall a,b\in \mathbb{R},\mathbb{\ }\mathfrak{n}(\left[ a,b\right) _{%
\mathbb{R}})=\left[ \left( b-a\right) +\varepsilon \right] \cdot \mathfrak{n}%
(\left[ 0,1\right) _{\mathbb{R}}  \label{due}
\end{equation}%
where $\varepsilon $ is an infinitesimal which might depend on $a$ and $b$.

Now we will define $\Lambda _{\mathbb{R}}\subset P_{fin}(\mathbb{R})$ in
such a way that eq.~(\ref{uno}) and eq.~(\ref{due}) be satisfied; first, we
set 
\begin{equation*}
\Theta =P_{fin}(\left[ 0,1\right] _{\mathbb{R}}\backslash \left[ 0,1\right]
_{\mathbb{Q}})
\end{equation*}%
namely $\Theta $ is the family of the finite sets of irrational numbers
between 0 and 1. Then for any $n\in \mathbb{N}$ and $\theta \in \Theta $, we
set 
\begin{equation*}
\mu _{n,\theta }:=\mu _{n}\cup \left\{ \frac{p+a}{n}\ |\ \frac{p}{n}\in \mu
_{n}\backslash \left\{ n\right\} ,\ a\in \theta \right\}
\end{equation*}
where $\mu _{n}$ is defined by eq.~(\ref{ln}). You may think to have
``constructed'' $\mu _{n,\theta }$ in the following way: you start with the
segment $\left[ -n,n\right] $ and you divide it in $n^{2}$ parts of length $%
1/n$ of the form $\left[ \frac{p}{n},\frac{p+1}{n}\right] ,$ $%
p=-n^{2},-n^{2}+1,....,n^{2}-1.$ In each of these parts you put a
``rescaled'' copy of $\theta ,$ namely points of the form $\frac{p+a}{n}$
with $a\in \theta .$\ Thus, the set $\mu _{n,\theta }$ contains $n^{2}+1\ $%
rational numbers and $n^{2}\cdot \left\vert \theta \right\vert $ irrational
numbers.

Then we set 
\begin{equation}
\Lambda _{\mathbb{R}}=\left\{ \mu _{n,\theta }\ |\ n=m!,\ m\in \mathbb{N},\
\theta \in \Theta \right\}.  \label{lambdaR}
\end{equation}

\begin{proposition}
If $\Lambda =\Lambda _{\mathbb{R}}$ is given by eq.~(\ref{lambdaR}), then
eq.~(\ref{uno}) and eq.~(\ref{due}) hold.
\end{proposition}

\textbf{Proof}. Take an interval $\left[ a,b\right) $ with $a,b\in \mathbb{Q}
$; then, if $n$ is sufficiently large $\left[ a,b\right) \cap \mu _{n,\theta
}$ contains $n\left( b-a\right) $ rational numbers and $n\left\vert \theta
\right\vert \left( b-a\right) $ irrational numbers so that 
\begin{equation*}
\left\vert \left[ a,b\right) \cap \mu _{n,\theta }\right\vert =n\left(
b-a\right) \left( \left\vert \theta \right\vert +1\right) .
\end{equation*}%
Then, if we choose $a=0$ and $b=1$, we have that 
\begin{equation*}
\mathfrak{n}(\left[ 0,1\right) _{\mathbb{R}})=\lim_{\mu _{n,\theta }\in
\Lambda _{\mathbb{R}}}\left\vert \left[ 0,1\right) \cap \mu _{n,\theta
}\right\vert =\lim_{\mu _{n,\theta }\in \Lambda _{\mathbb{R}}}n\left(
\left\vert \theta \right\vert +1\right) .
\end{equation*}%
Then, in general we have that 
\begin{eqnarray*}
\mathbb{\ }\mathfrak{n}(\left[ a,b\right) _{\mathbb{R}}) &=&\lim_{\mu
_{n,\theta }\in \Lambda _{\mathbb{R}}}\left\vert \left[ a,b\right) \cap \mu
_{n,\theta }\right\vert =\lim_{\mu _{n,\theta }\in \Lambda _{\mathbb{R}}}%
\left[ n\left( b-a\right) \left( \left\vert \theta \right\vert +1\right) %
\right] \\
&=&\left( b-a\right) \cdot \lim_{\mu _{n,\theta }\in \Lambda _{\mathbb{R}%
}}n\left( \left\vert \theta \right\vert +1\right) =\left( b-a\right) 
\mathfrak{n}(\left[ 0,1\right) _{\mathbb{R}}).
\end{eqnarray*}%
Then eq.~(\ref{uno}) holds. Now let us prove eq.~(\ref{due}).

If $a\in \mathbb{R}\setminus \mathbb{Q\ }$or $b\in \mathbb{R}\setminus 
\mathbb{Q}$, then $\left[ a,b\right) \cap \mu _{n}$ contains at most $%
n\left( b-a\right) +1\ $rational numbers, in fact, in $\left[ a,b\right) \ $%
you can fit at most $n\left( b-a\right) $ intervals of the form $\left[ 
\frac{p}{n},\frac{p+1}{n}\right] .$ Moreover, since you have at most $%
n\left( b-a\right) $ intervals of the form $\left[ \frac{p}{n},\frac{p+1}{n}%
\right] $ and two smaller intervals at the extremes of the form $\left[ a,%
\frac{p_{a}}{n}\right] $ and $\left[ \frac{p_{b}}{n},b\right] $ for a
suitable choice of $p_{a}$ and $p_{b},$ $\left[ a,b\right) \cap \mu
_{n,\theta }\ $contains at most $\left[ n\left( b-a\right) +2\right]
\left\vert \theta \right\vert $ irrational numbers; then 
\begin{eqnarray*}
\left\vert \left[ a,b\right) \cap \mu _{n,\theta }\right\vert &\leq &n\left(
b-a\right) +1+\left[ n\left( b-a\right) +2\right] \left\vert \theta
\right\vert \\
&=&n\left( b-a\right) \left( \left\vert \theta \right\vert +1\right)
+2\left\vert \theta \right\vert +1 \\
&\leq &n\left( b-a+\frac{2}{n}\right) \left( \left\vert \theta \right\vert
+1\right) .
\end{eqnarray*}%
Thus 
\begin{eqnarray*}
\mathbb{\ }\mathfrak{n}(\left[ a,b\right) _{\mathbb{R}}) &=&\lim_{\mu
_{n,\theta }\in \Lambda _{\mathbb{R}}}\left\vert \left[ a,b\right) \cap \mu
_{n,\theta }\right\vert \leq \lim_{\mu _{n,\theta }\in \Lambda _{\mathbb{R}}}%
\left[ n\left( b-a+\frac{2}{n}\right) \left( \left\vert \theta \right\vert
+1\right) \right] \\
&\leq &\lim_{\mu _{n,\theta }\in \Lambda _{\mathbb{R}}}\left( b-a+\frac{2}{n}%
\right) \cdot \lim_{\mu _{n,\theta }\in \Lambda _{\mathbb{R}}}n\left(
\left\vert \theta \right\vert +1\right) \\
&=&\left( b-a+\frac{2}{\alpha }\right) \mathfrak{n}(\left[ 0,1\right) _{%
\mathbb{R}}).
\end{eqnarray*}

Moreover, since $\left[ a,b\right) $ contains at least $n\left( b-a\right)
-1 $ interval of type $\left[ \frac{p}{n},\frac{p+1}{n}\right] ,\ $arguing
in a similar way as before, we have that 
\begin{eqnarray*}
\left\vert \left[ a,b\right) \cap \mu _{n,\theta }\right\vert &\geq &n\left(
b-a\right) +\left[ n\left( b-a\right) -1\right] \left\vert \theta \right\vert
\\
&=&n\left( b-a\right) \left( \left\vert \theta \right\vert +1\right)
-\left\vert \theta \right\vert \\
&\geq &n\left( b-a-\frac{1}{n}\right) \left( \left\vert \theta \right\vert
+1\right).
\end{eqnarray*}
Then, 
\begin{eqnarray*}
\mathbb{\ }\mathfrak{n}(\left[ a,b\right) _{\mathbb{R}}) &=&\lim_{\mu
_{n,\theta }\uparrow \mathbb{R}}\left\vert \left[ a,b\right) \cap \mu
_{n,\theta }\right\vert \geq \lim_{\mu _{n,\theta }\uparrow \mathbb{R}}\left[
n\left( b-a-\frac{1}{n}\right) \left( \left\vert \theta \right\vert
+1\right) \right] \\
&\geq &\lim_{\mu _{n,\theta }\uparrow \mathbb{R}}\left( b-a-\frac{1}{n}%
\right) \cdot \lim_{\mu _{n,\theta }\uparrow \mathbb{R}}n\left( \left\vert
\theta \right\vert +1\right) \\
&=&\left( b-a-\frac{1}{\alpha }\right) \mathfrak{n}(\left[ 0,1\right) _{%
\mathbb{R}}).
\end{eqnarray*}
Thus eq.~(\ref{due}) follows with 
\begin{equation*}
\left\vert \varepsilon \right\vert \leq \frac{2}{\alpha }.
\end{equation*}

$\square $

So, if we have two intervals $\left[ a_{0},b_{0}\right) _{\mathbb{R}}\subset %
\left[ a_{1},b_{1}\right) _{\mathbb{R}},$ using the above proposition, we
get \ref{tre}. Moreover, it is easy to prove that the NAP-space produced by $%
\left( \mathbb{R},1,I_{\Lambda _{\mathbb{R}}}\right) $ is consistent with $%
\left( \mathbb{Q},1,I_{\Lambda _{\mathbb{Q}}}\right) $, namely that the
analog of eq.~(\ref{PNA}) holds: if $A\subset \mathbb{Q}$, then 
\begin{equation*}
P_{\mathbb{Q}}(A)=P_{\mathbb{R}}(A\ |\ \mathbb{Q}).
\end{equation*}

\subsection{The infinite sequence of coin tosses}

Let us consider an infinite sequence of tosses with a fair coin.\footnote{%
This example is used by Williamson in an attempt to refute the possibility
of assigning infinitesimal probability values to a particular outcome of
such a sequence \cite{Williamson2007}. As observed by Weintraub,
Williamson's argument relies on a relabeling of the individual tosses, which
is not compatible with non-Archimedean probabilities \cite{Weintraub2008}.}
In the Kolmogorovian framework, the infinite sequence of fair coin tosses is
modeled by the triple $\left( \Omega ,\mathfrak{A},\mu \right)$, where $%
\Omega =\left\{H,T\right\} ^{\mathbb{N}}$ is the space of sequences which
take values in the set $\left\{ H,T\right\} $ namely \textit{Heads} and 
\textit{Tails}. We will denote by $\omega =\left( \omega _{1},...,\omega
_{n},.....\right) $ the generic sequence.

$\mathfrak{A}$ is the $\sigma $-algebra generated by the `cylindrical sets'.
A cylindrical set of codimension $n$ is defined by a $n$-ple of indices $%
\left( i_{1},...,i_{n}\right) $ and an $n$-ple of elements in $%
\left\{H,T\right\} ,$ namely $\left( t_{1},...,t_{n}\right) $ where $t_{k}\ $%
is either $H$ or $T.$

A cylindrical set of codimension $n$ is defined as follows: 
\begin{equation*}
C_{\left( t_{1},...,t_{n}\right) }^{\left( i_{1},...,i_{n}\right) }=\left\{
\omega \in \Omega \ |\ \omega _{i_{k}}=t_{k}\right\}.
\end{equation*}

From the probabilistic point of view, $C_{\left(
t_{1},...,t_{n}\right)}^{\left( i_{1},...,i_{n}\right) }$ represents the
event that that $i_{k}$-th coin toss gives $t_{k}\ $for $k=1,...,n.$

The probability measure on the generic cylindrical set is given by 
\begin{equation}
\mu \left( C_{\left( t_{1},...,t_{n}\right) }^{\left( i_{1},...,i_{n}\right)
}\right) =2^{-n}.  \label{pcil}
\end{equation}%
The measure $\mu$ can be extended in a unique way to $P_{KA}$ on the algebra 
$\mathfrak{A}$ (by Caratheodory's theorem).

In this particular model, you can see the problems with the Kolmogorovian
approach which we discussed in section~\ref{sec:PK}:

\begin{itemize}
\item every event $\left\{ \omega \right\} \in \Omega $ has $0$ probability
but the union of all these `seemingly impossible' events has probability $1$;

\item if $F$ is a finite set, and $\left\{ \omega \right\} \subset F,$ then
the conditional probability $P_{KA} (\left\{ \omega \right\} |F)$ is not
defined; nevertheless, you know that the conditionalizing event is not the
empty event, so the conditional probability should be defined: it makes
sense and its value is $\frac{1}{\left\vert F\right\vert }$;

\item there are subsets of $\Omega$ for which the probability is not defined
(namely the non-measurable sets).
\end{itemize}

Now, we will construct a NAP so that we can compare the two different
approaches. We need to construct a NAP produced by $\left( \Omega
,w,I_{\Lambda} \right)$ which satisfies the following assumptions:

\begin{itemize}
\item (i) if $F\subset \Omega $ is a finite non-empty set, then 
\begin{equation*}
P(A|F)=\frac{\left\vert A\cap F\right\vert }{\left\vert F\right\vert };
\end{equation*}

\item (ii) eq.~(\ref{pcil}) holds for $P,$ namely 
\begin{equation}
P\left( C_{\left( t_{1},...,t_{n}\right) }^{\left( i_{1},...,i_{n}\right)
}\right) =2^{-n}.  \label{napcil}
\end{equation}
\end{itemize}

\bigskip

Experimentally, we can only observe a finite numbers of outcomes: both
cylindrical events (\textit{cf.}~(ii)) and finite conditional probability (%
\textit{cf.}~(i)) are based on a finite number of observations. In some
sense, (i) and (ii) are the `experimental data' on which to construct the
model.

Property (i) implies that we get a fair probability, thus we have to take $%
w\equiv 1;$ so every infinite sequence of coin tosses in $\Omega
=\left\{H,T\right\} ^{\mathbb{N}}$ has probability $1/\mathfrak{n}(\Omega )$%
. Property (ii) is the analog of eq.~(\ref{uno}) in the case of a fair
lottery on $\mathbb{R}$.\footnote{%
Eq.~(\ref{napcil}) implies that for any $\mu $-measurable set $E$, we have
that $P\left( E\right) \sim \mu (E)$ and this relation is the analog of eq.~(%
\ref{tre}).} So, we have to choose $\Lambda $ in such a way that eq.~(\ref%
{napcil}) holds.

\bigskip

To do this we need some other notation; if $b=\left( b_{1},...,b_{n}\right)
\in \left\{ H,T\right\} ^{n}$ is a finite string and $c=\left(
c_{1},...,c_{n},.....\right) \in \left\{ H,T\right\} ^{\mathbb{N}}$ is an
infinite sequence, then we set 
\begin{equation*}
b\circledast c=\left( b_{1},...,b_{n},c_{1},...,c_{k},.....\right)
\end{equation*}
namely, the sequence $b\mathfrak{\circledast }c$ is obtained by the sequence 
$b$ followed by the infinite sequence $c$. Now, if $\sigma \in \mathcal{P}%
_{fin}\left( \left\{ H,T\right\} ^{\mathbb{N}}\right) $ and $n\in \mathbb{N}$%
, we set 
\begin{equation*}
\lambda _{n,\sigma }=\left\{ b\mathfrak{\circledast }c\ |\ b\in
\left\{H,T\right\} ^{n}\ \text{and}\ c\in \sigma \right\}
\end{equation*}
and 
\begin{equation*}
\Lambda _{CT}:=\left\{ \lambda _{n,\sigma }\ |\ \sigma \in \mathcal{P}%
_{fin}\left( \left\{ H,T\right\} ^{\mathbb{N}}\right) \ \text{and\ }n\in 
\mathbb{N}\right\} .
\end{equation*}
Notice that $\left( \left\{ H,T\right\} ^{\mathbb{N}},1,I_{\Lambda_{CT}}%
\right)$ produces a well-defined NAP-space since $\Lambda_{CT}$ is a
directed set; in fact 
\begin{equation*}
\lambda _{n_{1},\sigma _{1}}\cup \lambda _{n_{2},\sigma _{2}}\subset \lambda
_{\max (n_{1},n_{2}),\sigma _{3}}
\end{equation*}
for a suitable choice of $\sigma _{3}\in \mathcal{P}_{fin}\left(
\left\{H,T\right\} ^{\mathbb{N}}\right) .$

Moreover $\Lambda _{CT}$ is the `wise choice' as the following theorem shows:

\begin{theorem}
If $P$ is the NAP produced by $\left( \left\{ H,T\right\} ^{\mathbb{N}%
},1,I_{\Lambda_{CT}}\right) $, then eq.~(\ref{napcil}) holds.
\end{theorem}

\textbf{Proof. }We have 
\begin{equation}
\mathfrak{n}(\Omega )=\lim_{\lambda _{N,\sigma }\in \Lambda _{CT}}\left\vert
\lambda _{N,\sigma }\right\vert =\lim_{\lambda _{N,\sigma }\in \Lambda
_{CT}}\left( 2^{N}\cdot \left\vert \sigma \right\vert \right).  \label{no}
\end{equation}

Now consider the cylinder $C_{\left( t_{1},...,t_{n}\right) }^{\left(
i_{1},...,i_{n}\right) }$ and take $N=\max \ i_{n}.$ Then, for every $\sigma$%
, we have that 
\begin{equation*}
\lambda _{N,\sigma }\cap C_{\left( t_{1},...,t_{n}\right) }^{\left(
i_{1},...,i_{n}\right) }=\left\{ \omega \in \lambda _{N,\sigma }\ |\ \omega
_{i_{k}}=t_{k}\right\}.
\end{equation*}
Then 
\begin{equation*}
\left\vert \lambda _{N,\sigma }\cap C_{\left( t_{1},...,t_{n}\right)
}^{\left( i_{1},...,i_{n}\right) }\right\vert =\frac{\left\vert \lambda
_{N,\sigma }\right\vert }{2^{n}}=2^{N-n}\cdot \left\vert \sigma \right\vert
\end{equation*}
and hence, by eq.~(\ref{no}), 
\begin{eqnarray*}
\mathfrak{n}\left( C_{\left( t_{1},...,t_{n}\right)
}^{\left(i_{1},...,i_{n}\right) }\right) &=&\lim_{\lambda _{N,\sigma }\in
\Lambda _{CT}}\left\vert \lambda _{N,\sigma }\cap C_{\left(
t_{1},...,t_{n}\right) }^{\left( i_{1},...,i_{n}\right) }\right\vert
=\lim_{\lambda _{N,\sigma }\in \Lambda _{CT}}\left( 2^{N-n}\cdot \left\vert
\sigma \right\vert \right) \\
&=&2^{-n}\lim_{\lambda _{N,\sigma }\in \Lambda _{CT}}\left( 2^{N}\cdot
\left\vert \sigma \right\vert \right) =2^{-n}\cdot \mathfrak{n}(\Omega).
\end{eqnarray*}
Concluding, we have that 
\begin{equation*}
P\left( C_{\left( t_{1},...,t_{n}\right) }^{\left( i_{1},...,i_{n}\right)
}\right) =\frac{\mathfrak{n}\left( C_{\left( t_{1},...,t_{n}\right)
}^{\left( i_{1},...,i_{n}\right) }\right) }{\mathfrak{n}(\Omega )}=2^{-n}.
\end{equation*}

$\square $

\end{document}